\documentclass[aap,MSNbibl,nameyear,seceqn,dvips]{arximspdf}
\usepackage{algorithmic}
\usepackage{algorithm}

\doi{10.1214/13-AAP974} %kopijuoti is PTS
\volume{24}
\issue{5}
\pubyear{2014}
\firstpage{2143}
\lastpage{2175}

\makeatletter
\setattribute{abstract}{width}{292pt}
\newcommand{\widebar}{\overline}
\newcommand{\rrvert}{\vert}
\newcommand{\llvert}{\vert}
\newcommand{\eqref}[1]{(\ref{#1})}
\newcommand{\bfnu}{{\bolds\nu}}
\newcommand{\integer}{{\mathbb Z}}
\newcommand{\pintegers}{{\integer_+}}
\def\bP{\mathbf{P}}
\def\argmin{\arg\min}
\newtheorem{Th}{Theorem}[section]
\newtheorem{Lm}{Lemma}[section]
\newproclaim{Rk}{Remark}[section]
\newproclaim{hy}{Hypotheses}
\makeatother

\begin{document}
\begin{frontmatter}

\title{Rare event simulation for processes generated via stochastic
fixed point equations}
\runtitle{Importance sampling for SFPE}

\begin{aug}
\author[A]{\fnms{Jeffrey F.} \snm{Collamore}\corref{}\thanksref{t1}\ead[label=e1]{collamore@math.ku.dk}},
\author[B]{\fnms{Guoqing} \snm{Diao}\ead[label=e2]{gdiao@gmu.edu}}\\
\and
\author[B]{\fnms{Anand N.} \snm{Vidyashankar}\thanksref{t2}\ead[label=e3]{avidyash@gmu.edu}}
\runauthor{J. F. Collamore, G. Diao and A. N. Vidyashankar}
\affiliation{University of Copenhagen, George Mason University\\ and
George Mason University}
\address[A]{J. F. Collamore\\
Department of Mathematical Sciences\\
University of Copenhagen\\
Universitetsparken 5\\
DK-2100 Copenhagen {\O}\\
Denmark\\
\printead{e1}} %adresu isvedimo komanda gale!
\address[B]{G. Diao\\
A. N. Vidyashankar\\
Department of Statistics\\
George Mason University\\
4400 University Drive,
MS 4A7\\
Fairfax, Virginia 22030\\
USA\\
\printead{e2}\\
\phantom{E-mail:\ }\printead*{e3}}
\end{aug}
\thankstext{t1}{Supported in part by Danish Research Council (SNF)
Grant ``Point Process Modelling and Statistical Inference,'' No. 09-092331.}
\thankstext{t2}{Supported by Grant NSF DMS 000-03-07057.}

% HISTORY:
\received{\smonth{7} \syear{2011}}
\revised{\smonth{9} \syear{2013}}

% ABSTRACT
%
\begin{abstract}
In a number of applications, particularly in financial and actuarial
mathematics, it is of interest to characterize
the tail distribution of a random variable $V$ satisfying the
distributional equation $V \stackrel{\mathcal D}{=} f(V)$,
where $f(v) = A \max\{v,D\} + B$ for $(A,B,D) \in(0,\infty) \times
{\mathbb R}^2$.
This paper is concerned with computational methods for evaluating these
tail probabilities.
We introduce a novel importance sampling algorithm, involving an
exponential shift over
a random time interval, for estimating these rare event
probabilities. We prove that the proposed estimator is: (i)~consistent,
(ii)~strongly efficient and (iii) optimal within
a wide class of dynamic importance sampling estimators. Moreover, using
extensions of ideas from nonlinear renewal theory,
we provide a precise description of the running time of the algorithm.
To establish these results, we develop new
techniques concerning the convergence of moments of stopped perpetuity
sequences, and the first entrance
and last exit times of associated Markov chains on ${\mathbb R}$.
We illustrate our methods with a variety of numerical examples which
demonstrate the ease and scope of the implementation.
\end{abstract}

% KEYWORDS
% Pirmas kwd is didziosios raides
%

\end{frontmatter}

%s1 #&#
\section{Introduction}\label{sec1}
%%renewcommand{\baselinestretch}{1.5}
This paper introduces a rare event simulation algorithm for
estimating the tail probabilities of the \textit{stochastic fixed point
equation} (SFPE)
%
%e1.1 #&#
\begin{equation}
\label{intro1} V \stackrel{{\mathcal D}} {=} f(V)\qquad\mbox{where }f(v)
\equiv A \max\{v,D\} + B
\end{equation}
for $(A,B,D) \in(0,\infty) \times{\mathbb R}^2$.
SFPEs of this general form arise in a wide variety of applications,
such as extremal estimates for financial time series models
and ruin estimates in actuarial mathematics.
%such as the ruin problem with investments and the valuation of large
%losses for pension funds.
Other related applications arise in branching processes in random environments
and the study of algorithms in computer science. See \citet{JC09},
\citet{JCAV11}, or Section~\ref{sec4} below for a more
detailed description of some of these applications.

In a series of papers [e.g., \citet{HK73}, \citet{WV79},
\citet{CG91}], the tail probabilities for the
SFPE~\eqref{intro1} have been asymptotically characterized. Under
appropriate moment and regularity conditions, it is known that
%
%e1.2 #&#
\begin{equation}
\label{intro3} \lim_{u \rightarrow\infty} u^{\xi}{\mathbf P} \{ V>u \} =C
\end{equation}
for finite positive constants $C$ and $\xi$, where $\xi$ is identified
as the nonzero solution
to the equation ${\mathbf E} [A^\alpha] = 1$.
Recently, in \citet{JCAV11}, the constant $C$ has been identified
as the
$\xi$th moment of the difference of a perpetuity sequence
and a conjugate sequence.

The purpose of this article is to introduce a rigorous computational
approach, based on importance sampling,
for Monte Carlo estimation of the rare event probability ${\mathbf P} \{ V
> u \}$. While importance sampling methods have been
developed for numerous large deviation problems involving i.i.d. and
Markov-dependent
random walks [cf. \citet{SAPG07}],
the adaptation of these methods to~\eqref{intro1} is distinct and
requires new techniques.
In this paper, we propose a nonstandard approach
involving a \textit{dual} change of measure of a process $\{ V_n\}$
performed over two random time intervals: namely, the excursion of $\{
V_n \}$
to $(u,\infty)$ followed by the return of this process to a given set
${\mathcal C} \subset{\mathbb R}$.

The motivation for our algorithm stems from the observation that the
SFPE~\eqref{intro1} induces a forward recursive sequence, namely,
%
%e1.3 #&#
\begin{equation}
\label{eq07} V_n = A_n \max\{ D_n,
V_{n-1} \} + B_n,\qquad n=1,2,\ldots, V_0=v,
\end{equation}
where $\{ (A_n,B_n,D_n)\dvtx  n \in\pintegers\}$ is an i.i.d. sequence
with the same law as $(A,B,D)$.
It is important to observe that in many applications,
the mathematical process under study is obtained through the \textit{backward} iterates of the given SFPE [as described by
\citet{GL86} or \citet{JCAV11}, Section~2.1]. For example,
the linear recursion $f(v) = Av + B$ induces the backward recursive sequence
or \textit{perpetuity sequence}
%
%e1.4 #&#
\begin{equation}
\label{eq08} Z_n:= V_0 + \frac{B_1}{A_1} +
\frac{B_2}{A_1 A_2} + \cdots+ \frac{B_n}{A_1 \cdots A_n},\qquad
n=1,2,\ldots.
\end{equation}
However, since $\{ Z_n \}$ is not Markovian, it is less natural to simulate
$\{ Z_n \}$ than the corresponding forward sequence $\{ V_n \}$. Thus,
a central aspect of our
approach is the \textit{conversion} of the given perpetuity sequence,
via its SFPE, into a forward recursive sequence which we then simulate.
Because $\{ V_n \}$ is Markovian, we can then study this process over
excursions emanating from, and then returning to,
a given set ${\mathcal C} \subset{\mathbb R}$.

In the special case of the perpetuity sequence in~\eqref{eq08},
simulation methods for estimating
${\mathbf P} \{ \lim_{n \to\infty} Z_n > u \}$ have recently been studied
in \citet{JBHLBZ11} under the strong
assumption that $\{ B_n \}$ is nonnegative. Their method is very
different from ours,
involving the simulation of $\{ Z_n \}$ \textit{directly} until the first\vadjust{\goodbreak}
passage time to a level $cu$, where $c \in(0,1)$, and
a rough analytical approximation to relate this probability to the
first passage probability at level $u$.
Their methods do not generalize to the other processes studied in this
paper, such as the ruin problem with investments
or related extensions.
In contrast, our goal here is to develop a \textit{general} algorithm
which is flexible and can be applied to the wider class of processes
governed by~\eqref{intro1} and some of its extensions.
While we focus on~\eqref{intro1}, it is worthwhile to mention here that
our algorithm provides an important ingredient for addressing a larger
class of problems, including nonhomogeneous recursions on trees, which
are analyzed in \citet{JCAV13a}. Also, it seems
plausible that the method should extend to the class of random maps
which can
be approximated by~\eqref{intro1} in the sense of \citet{JCAV11},
Section~2.4. This extension would encompass several other problems of applied
interest, such as the AR(1) process with ARCH(1) errors. Yet another
feasible generalization is to Markov-dependent recursions under Harris
recurrence, utilizing the reduction to i.i.d. recursions described
in \citet{JC09} and \citet{JCAV13}, Section~3.

In this paper, we present an algorithm and establish that it is
consistent and efficient; that is, it displays the bounded relative
error property.
It is interesting to note that in the proof of efficiency, certain new
issues arise concerning the convergence of the perpetuity sequence
\eqref{eq08}.
Specifically, while it is known that~\eqref{eq08} converges to a finite
limit under minimal conditions, the necessary and sufficient condition for
the $L_\beta$ convergence of $\{ Z_n \}$ in~\eqref{eq08} is that ${\mathbf
E} [ A^{-\beta} ] < 1$; cf. \citet{GAAIUR09}. However, our analysis
will involve moments of quantities similar to~$\{ Z_n \}$, but where
${\mathbf E} [ A^{-\beta} ]$ is \textit{greater} than one, and hence our perpetuity
sequences will necessarily be divergent in $L_\beta$.
To circumvent this difficulty, we
study these perpetuity sequences over randomly stopped intervals,
namely, over cycles
emanating from, and returning to, a given subset ${\mathcal C}$ of
${\mathbb R}$. As a technical point,
it is worth noting that if the return time, $K$, were replaced by the
more commonly studied regeneration time $\tau$ of the chain $\{ V_n \}
$, then the existing literature
on Markov chain theory would still not shed much light on the tails of
$\tau$ and hence the convergence of $V_\tau$.
Thus, the fact that $K$ has sufficient exponential tails for the
convergence of~$V_K$ is due to the recursive structure of the
particular class of Markov
chains we consider and seems to be a general property for this class of
Markov chains.
These results concerning the moments of $L_\beta$-divergent perpetuity
sequences complement the known literature on perpetuities and appear
to be of some independent interest.

Next, we go beyond the current literature by establishing a \textit{sharp}
asymptotic estimate for the running time of the
algorithm, thereby showing that our algorithm is, in fact, strongly
efficient; cf. Remark~\ref{re2.2} below. To this end, we introduce methods from
nonlinear renewal theory, as well as methods from Markov chain theory
involving the first entrance and last exit times of the process $\{ V_n
\}$.
Finally, motivated by the Wentzell--Freidlin theory of large
deviations, we provide an optimality result;
specifically, we consider other possible level-dependent changes of
measure for the process $\{ V_n \}$
selected from a wide class of dynamic importance sampling algorithms
[in the sense of \citet{PDHW05}].
We show that our algorithm is the \textit{unique} choice which attains
bounded relative error, thus establishing the validity of our
method amongst a natural class of possible algorithms.

%
% STATEMENT OF RESULTS
%
%s2 #&#
\section{The algorithm and a statement of the main results}\label{sec2}
%s2.1 #&#
\subsection{Background: The forward and backward recursive sequences}\label{sec2.1}

We start with a general SFPE of the form
%
%e2.1 #&#
\begin{equation}
\label{letac1} V \stackrel{\mathcal D} {=} f(V) \equiv F_Y(V),
\end{equation}
where $F_Y\dvtx  {\mathbb R} \times{\mathbb R}^d \to{\mathbb R}$ is
deterministic, measurable and continuous in its first component.
Let $v $ be an element of the range of $F_Y$, and let $\{ Y_n \}$ be an
i.i.d. sequence of r.v.'s such that $Y_n \stackrel{\mathcal D}{=} Y$
for all $n$.
Then the \textit{forward sequence} generated by the SFPE~\eqref{letac1} is
defined by
%
%e2.2 #&#
\begin{equation}
\label{letac2} V_n (v) = F_{Y_n} \circ F_{Y_{n-1}}
\circ\cdots\circ F_{Y_1}(v),\qquad n=1,2,\ldots,\ V_0 = v,
\end{equation}
whereas the \textit{backward sequence} generated by this SFPE is defined by
%
%e2.3 #&#
\begin{equation}
Z_n(v) = F_{Y_1} \circ F_{Y_2} \circ\cdots\circ F_{Y_n} (v),\qquad n=1,2,\ldots,\ Z_0 = v.
\end{equation}
While the forward sequence is always Markovian, the backward equation
need not be Markovian; however, for every
$v$ and $n$, $V_n(v)$ and $Z_n(v)$ are identically distributed. This
observation is critical since
it suggests that---regardless of whether the SFPE was originally
obtained via forward or backward iteration---a natural approach to analyzing
the process is through its forward iterates.

%s2.2 #&#
\subsection{Background: Asymptotic estimates}\label{sec2.2}
We now specialize to the recursion~\eqref{intro1}.
This recursion is often referred to as ``Letac's model~E.''\vadjust{\goodbreak}

Let $\mathfrak F_n$ denote the $\sigma$-field generated by $\{(A_i, B_i, D_i)\dvtx 1 \leq i \leq n\}$, and let
\[
\lambda(\alpha) = {\mathbf E} \bigl[ A^\alpha\bigr] \quad\mbox{and}\quad
\Lambda(\alpha) = \log\lambda(\alpha),\qquad\alpha\in{\mathbb R}.
\]
Let $\mu$ denote the distribution of $Y = (\log A,B,D)$ and $\mu_\alpha
$ denote the $\alpha$-shifted
distribution with respect to the first variable; that is,
%
%e2.4 #&#
\begin{equation}
\label{meas-shift} \mu_\alpha(E):= \frac{1}{\lambda(\alpha)} \int
_E e^{\alpha x}\,d\mu(x,y,z),\qquad E \in{\mathcal B}
\bigl({\mathbb R}^3\bigr), \alpha\in{\mathbb R},
\end{equation}
where, here and in the following, ${\mathcal B}(E)$ denotes the Borel
sets of $E$.
Let ${\mathbf E}_\alpha[ \cdot ]$ denote expectation with respect to this
$\alpha$-shifted measure.
% and for
%any $\beta$ in the domain of $\Lambda$, define

For any r.v. $X$, let
${\mathfrak L}(X)$ denote the probability law of $X$, and let $\operatorname{supp} (X)$ denote the support of $X$.
Also, write $X \sim{\mathfrak L}(X)$ to denote that
$X$ has this probability law. Given an i.i.d. sequence $\{ X_n \}$, we will
often write $X$ for a ``generic'' element of this sequence. Finally,
for any function $f$, let $\operatorname{dom} (f)$ denote the domain of $f$,
and let $f^\prime$, $f^{\prime\prime}$, etc. denote the
successive derivatives of $f$.
%We are now prepared to state the main hypotheses of the paper.

We now state the main hypotheses needed to establish the asymptotic
decay of ${\mathbf P} \{ V > u \}$ in~\eqref{intro3}; note that
\textup{(H$_0$)} is only needed to obtain the explicit\vadjust{\goodbreak} representation of $C$, as
given in \citet{JCAV11}. These conditions
will form the starting point of our study.

\begin{hy*} %\textit{Hypotheses}:
\begin{longlist}[\textup{(H$_3$)}]
\item[\textup{(H$_0$)}] The r.v. $A$ has an absolutely continuous component with
respect to Lebesgue measure with a nontrivial continuous density in
a neighborhood of ${\mathbb R}$.

\item[\textup{(H$_1$)}] $\Lambda(\xi) = 0$ for some $\xi\in(0,\infty)\cap \operatorname{dom} (\Lambda')$.

\item[\textup{(H$_2$)}] ${\mathbf E} [ | B |^\xi ] < \infty$ and ${\mathbf E} [ (A|D| )^\xi ] < \infty$.

\item[\textup{(H$_3$)}] ${\mathbf P} \{ A> 1, B> 0 \} > 0$
or ${\mathbf P} \{ A>1, B \ge0, D > 0 \} >0$.
\end{longlist}
\end{hy*}

Note that \textup{(H$_3$)} implies that the process $\{ V_n \}$ is nondegenerate
(i.e., it is not concentrated at a single point).

Under these hypotheses, it can be shown that the forward sequence $\{
V_n \}$ generated by the SFPE~\eqref{intro1}
is a Markov chain which is $\varphi$-irreducible and geometrically
ergodic [\citet{JCAV11}, Lemma 5.1].
Thus $\{ V_n \}$ \mbox{converges} to a r.v. $V$ \textit{which itself satisfies
the SFPE}~\eqref{intro1}.
Moreover, with respect to its $\alpha$-shifted measure, the process $\{
V_n \}$ is \textit{transient} [\citet{JCAV11}, Lemma 5.2].

Our present goal is to develop an efficient Monte Carlo algorithm for
evaluating ${\mathbf P} \{ V > u \}$, for \textit{fixed} $u$, which remains efficient
in the asymptotic limit
as $u \to\infty$.

%s2.3 #&#
\subsection{The algorithm}\label{sec2.3}
Since the forward process
$V_n = A_n \max\{ D_n,V_{n-1} \} + B_n$
satisfies $V_n \approx A_n V_{n-1}$ for large $V_{n-1}$, and since $\{
V_n \}$ is transient in its $\xi$-shifted measure,\vadjust{\goodbreak}
large deviation theory suggests that we consider shifted distributions
and, in particular, the shifted measure $\mu_\xi$, where $\xi$ is
given as in \textup{(H$_1$)}.
To relate \mbox{${\mathbf P} \{ V > u \}$} under its original measure
to the paths of $\{ V_n \}$ under $\mu_\xi$-measure, let ${\mathcal C}:
= [-M,M]$ for some $M \ge0$,
and let $\pi$ denote the stationary distribution of~$\{ V_n \}$. Now
define a probability measure $\gamma$ on ${\mathcal C}$ by setting
%
%e2.5 #&#
\begin{equation}
{\label{defSD}} \gamma(E) = \frac{\pi(E)}{\pi({\mathcal C})},\qquad E
\in{\mathcal B}({\mathcal
C}).
\end{equation}
Let $K:= \inf \{ n\in\pintegers\dvtx  V_n \in{\mathcal C} \}$. Then in
Section~\ref{sec3}, we will establish the following representation formula:
%
%e2.6 #&#
\begin{equation}
\label{excursion-motiv} {\mathbf P} \{ V > u \} = \pi({\mathcal C}) {\mathbf
E}_\gamma[
N_u ],\qquad N_u:= \sum
_{n=0}^{K-1} {\mathbf1}_{\{ V_n > u \}},
\end{equation}
where ${\mathbf E}_{\gamma}[\cdot]$ denotes the expectation when the
initial state $V_0 \sim\gamma$. Thus motivated by large deviation
theory and the previous formula,
we simulate $\{ V_n \}$ \textit{over a cycle} emanating from the set
${\mathcal C}$ (with initial state $V_0 \sim\gamma$), and
then returning to ${\mathcal C}$, where simulation is performed
in the dual measure, which we now describe.

Set $T_u = \inf\{ n\dvtx  V_n > u \}$, and let
{\renewcommand{\theequation}{${\mathfrak D}$}
%e2.7 #&#
\begin{equation}
\label{DUAL-MEAS} {\mathfrak L} (\log A_n,B_n,D_n
) = \cases{ \mu_\xi, &\quad for $n=1,\ldots,T_u$,
\cr
\mu, &\quad for $n > T_u$,}
\end{equation}}\setcounter{equation}{6}%
where $\mu_\xi$ is defined as in~\eqref{meas-shift} and $\xi$ is given
as in \textup{(H$_1$)}. Let $\{ V_n \}$ be generated by the forward recursion
\eqref{eq07}, but with a \textit{driving sequence} $\{ V_n \} \equiv\{
(\log A_n,B_n,D_n) \}$ which is governed by \eqref{DUAL-MEAS} rather
than by
the fixed measure $\mu$.
Roughly speaking, the ``dual measure'' \eqref{DUAL-MEAS} shifts the
distribution of $\log A_n$ on a path of $\{ V_n \}$ until this process
exceeds the level $u$,
and reverts to the original measure thereafter. Let ${\mathbf E}_{\mathfrak
D} [ \cdot ]$ denote expectation with respect to \eqref{DUAL-MEAS}.

To relate the simulated sequence in the dual measure
to the required probability in the original measure, we introduce a
weighting factor.
Specifically, in the proof of Theorem~\ref{MainTh-1} below, we will show
\[
{\mathbf E}_{\mathfrak D} [{\mathcal E}_u ]= \pi({\mathcal C}) {\mathbf
E}_{{\mathfrak D}} \bigl[ N_u e^{-\xi S_{T_u}} {\mathbf
1}_{\{ T_u < K \}} | V_0 \sim\gamma\bigr],
\]
where $S_n:= \sum_{i=1}^n \log A_i$ and $\gamma$ is given as in
\eqref{defSD}. Using this identity, it is natural
to introduce the \textit{importance sampling estimator}
%
%e2.7 #&#
\begin{equation}
\label{alg-10} {\mathcal E}_u = N_u e^{-\xi S_{T_u}}
{\mathbf1}_{\{ T_u < K\}}.
\end{equation}
Then $\pi({\mathcal C}) {\mathcal E}_u$ is an unbiased estimator for
${\mathbf P} \{ V>u \}$.
However, since the stationary distribution $\pi$ and hence the distribution
$\gamma$ is seldom known---even if the underlying distribution of
$(\log A,B,D)$ is known---we
first run multiple realizations of~$\{ V_n \}$ according to the known measure
$\mu$ and thereby estimate\vadjust{\goodbreak}
$\pi({\mathcal C})$ and $\gamma$. Let $\hat{\pi}_k({\mathcal C})$, $\hat
{\gamma}_k$ denote the estimates obtained for $\pi({\mathcal C})$,
$\gamma$,
respectively, and let $\widehat{\mathcal E}_{u,n}$ denote the estimate obtained
upon averaging the realizations of ${\mathcal E}_u$.
This yields the estimator~$\hat{\pi}_k({\mathcal C}) \widehat{\mathcal E}_{u,n}$.

This discussion can be formalized as follows:\vspace*{12pt}

\hrule\vspace{6pt}
\noindent Rare event simulation algorithm using forward iterations of the SFPE\vspace*{4pt}
\hrule\vspace{3pt}
\begin{algorithmic}
\STATE$V_0 \sim\hat{\gamma}_k, m=0$
\REPEAT
\STATE$m \leftarrow m+1$
\STATE$V_m = A_m \max\{D_m, V_{m-1}\}+B_m$, $(\log A_m, B_m, D_m)
\sim\mu_{\xi}$
\UNTIL{$V_m > u$ or $V_m \in{\mathcal C}$}
\IF{$V_m >u$}
\REPEAT
\STATE$m \leftarrow m+1$
\STATE$V_m = A_m \max\{D_m, V_{m-1}\}+B_m$, $(\log A_m, B_m, D_m)
\sim\mu$
\UNTIL{$V_m \in{\mathcal C}$}
\STATE${\mathcal E}_u = N_u e^{-\xi S_{T_u}} {\mathbf1}_{\{ T_u < K\}}$
\ELSE
\STATE${\mathcal E}_u=0$
\ENDIF\vspace{0.05in}
\hrule\vspace{6pt}
\end{algorithmic}

The\vspace*{1pt} actual estimate is then obtained by letting ${\mathcal E}_{u,j}$
$(j=1,\ldots,n)$ denote the realizations of ${\mathcal E}_u$ produced
by the algorithm and
setting ${\mathbf P} \{V > u \} = \hat{\pi}_k({\mathcal C}) \widehat{\mathcal E}_{u,n}$,
where
\[
\hat{\pi}_k({\mathcal C})= \frac{1}{k} \sum
_{j=1}^k {\mathbf1}_{\{V^{(j)} \in{\mathcal C}\}} \quad\mbox{and}\quad\widehat{\mathcal E}_{u,n} = \frac{1}{n} \sum
_{j=1}^n {\mathcal E}_{u,j},
\]
where $V^{(1)}, V^{(2)}, \ldots, V^{(k)}$ is a sample from the
distribution of $V$ (which, we emphasize, is sampled from the \textit{center} of the distribution). In Section~\ref{sec4}, we describe how to obtain
samples from $V$ from a practical perspective. Finally, note that $\widehat{\mathcal E}_{u,n}$ also depends on $k$.

It is worth observing that in the special case $D=1$ and $B=0$, Letac's
model~E reduces to a multiplicative random walk. Moreover, in that case,
one can always take $\gamma$ to be a point mass at $\{ 1 \}$, at which
point the process regenerates. In this much-simplified setting, our algorithm
reduces to a standard regenerative importance sampling algorithm, as
may be used to evaluate the stationary exceedance probabilities in a
GI/G/1 queue.

%s2.4 #&#
\subsection{Consistency and efficiency of the algorithm}\label{sec2.4}
We begin by stating our results on consistency and efficiency.

%
%th2.1 #&#
\begin{Th}\label{ConsistencyTh} Assume Letac's model~E, and
suppose that \textup{(H$_1$)}, \textup{(H$_2$)} and~\textup{(H$_3$)} are satisfied.\vadjust{\goodbreak} Then for any
${\mathcal C}$
such that ${\mathcal C} \cap \operatorname{supp} (\pi) \neq\varnothing$ and any
$u$ such that $u \notin{\mathcal C}$, the algorithm is strongly
consistent; that is,
%
%e2.8 #&#
\begin{equation}
\label{Cons} \lim_{ k \to\infty} \lim_{n \to\infty} \hat{
\pi}_{k} ({\mathcal C}) \widehat{\mathcal E}_{u,n} = {\mathbf P} \{
V > u \}\qquad\mbox{a.s.}
\end{equation}
\end{Th}

%
%re2.1 #&#
\begin{Rk}
If the stationary distribution $\pi$ of $\{ V_n \}$ is
known on ${\mathcal C}$ (e.g., ${\mathcal C}=\{v\}$ for $v \in{\mathbb
R}$), then it will follow from the proof
of the theorem that
$\pi({\mathcal C}) \widehat{\mathcal E}_{u,n}$ is an unbiased estimator for
${\mathbf P} \{ V>u \}$.
\end{Rk}

%
%th2.2 #&#
\begin{Th} \label{MainTh-1}
Assume Letac's model~E, and suppose that \textup{(H$_1$)} and \textup{(H$_3$)} are
satisfied. Also, in place of \textup{(H$_2$)}, assume that for some $\alpha> \xi$,
%
%e2.9 #&#
\begin{equation}
\label{newcond-thm22} {\mathbf E} \bigl[ \bigl( A^{-1} | B |^2
\bigr)^{\alpha} \bigr] < \infty\quad\mbox{and}\quad{\mathbf E} \bigl[
\bigl(A|D|^2 \bigr)^\alpha\bigr] < \infty.
\end{equation}
Moreover, assume that one of the following two conditions holds:
$\lambda(\alpha) < \infty$ for some
$\alpha< -\xi$; or ${\mathbf E} [ (|D| + (A^{-1}|B|)) ^\alpha ] < \infty$
for all $\alpha> 0$. Then,
there exists an $M >0$ such that
%
%e2.10 #&#
\begin{equation}
\label{BRE} \sup_{u \ge0} \sup_{k \in\pintegers}
u^{2\xi} {\mathbf E}_{\mathfrak D} \bigl[ {\mathcal E}_u^2
| V_0 \sim\hat{\gamma}_k \bigr] < \infty.
\end{equation}
\end{Th}

Equation~\eqref{BRE} implies that our estimator exhibits bounded
relative error.
However,
a good choice of $M$ is critical for the practical usefulness of the
algorithm. A canonical method for choosing $M$ can be based on the
drift condition satisfied
by $\{ V_n \}$ (as given in Lemma~\ref{le3.1} below), but in practice, a proper
choice of $M$ is problem-dependent and only obtained numerically based
on the methods we introduce below in Section~\ref{sec4}.

%
% Running Time of the algorithm

%s2.5 #&#
\subsection{Running time of the algorithm}\label{sec2.5}
Next we provide precise asymptotics for the running time of the algorithm.
In the following theorem, recall that $K$ denotes the first return time
to ${\mathcal C}$ (corresponding to the termination of the algorithm), whereas
$T_u$ denotes the first passage time to $(u, \infty)$.

%
%th2.3 #&#
\begin{Th} \label{MainTh-2}
Assume Letac's model~E, and suppose that hypotheses \mbox{\textup{(H$_0$)--(H$_3$)}}
hold, $\Lambda^{\prime\prime\prime}$ is finite on $\{ 0,\xi\}$ and for
some $\varepsilon> 0$,
%
%e2.11 #&#
\begin{equation}
\label{march22f} {\mathbf P}_{\xi} \{V_1 \le1
|V_0=v \}=o\bigl(v^{-\varepsilon}\bigr)\qquad\mbox{as } v \to
\infty.
\end{equation}
Then
%
%e2.12 #&#
%e2.13 #&#
\begin{eqnarray}
\label{thm23a} {\mathbf E}_{\mathfrak D} [ K {\mathbf1}_{\{ K < \infty\}} ] &
<& \infty;
\\
\label{thm23b} \lim_{u \to\infty} {\mathbf E}_{\mathfrak D}
\biggl[\frac{T_u }{\log u} \bigg| T_u < K \biggr] & =& \frac{1}{\Lambda^\prime
(\xi)};
\end{eqnarray}
%
%
%e2.14 #&#
\begin{eqnarray}
\label{thm23c} \lim_{u \to\infty} {\mathbf E}_{\mathfrak D} \biggl[
\frac{K-T_u }{\log u}\bigg| T_u < K \biggr] & =& \frac{1}{|\Lambda^\prime(0)|}.
\end{eqnarray}
\end{Th}

%
%re2.2 #&#
\begin{Rk}\label{re2.2}
The ultimate objective of the algorithm is to minimize the simulation
cost, that is, the total number of Monte Carlo simulations needed to
attain a given
accuracy.
This grows according to
%
%e2.15 #&#
\begin{equation}
\label{newcrit}\quad  \operatorname{Var} ({\mathcal E}_{u} ) \bigl
\{c_1{\mathbf E}_{\mathfrak D} [ K | T_u < K ] +
c_2{\mathbf E}_{\mathfrak D} [ K {\mathbf1}_{\{T_u \ge K \}} ] \bigr\}\qquad\mbox{as } u \to\infty
\end{equation}
for appropriate constants $c_1$ and $c_2$;
cf. \citet{DS76}.
However, as a consequence of Theorem~\ref{THoptimality}, we have that under the dual
measure \eqref{DUAL-MEAS},
\[
{\mathbf E}_{\mathfrak D} [ K | T_u < K ] \sim\Theta\log u \qquad\mbox{as } u \to\infty
\]
for some positive constant $\Theta$,
while the last term in~\eqref{newcrit} converges to a finite constant.
Thus, by combining Theorems \ref{MainTh-2} and \ref{THoptimality}, we conclude that our algorithm
is indeed strongly efficient.
\end{Rk}

%s2.6 #&#
\subsection{Optimality of the algorithm}\label{sec2.6}
We conclude with a comparison of our algorithm to other algorithms
obtained through forward iterations
involving alternative measure transformations. A natural alternative
would be to simulate with some measure $\mu_\alpha$ until the time $T_u
= \inf\{ n\dvtx  V_n > u \}$ and revert to some other measure $\mu_\beta$
thereafter. More generally, we may consider simulating from a general
class of distributions with some form of state dependence, as we now describe.

Let $\nu(\cdot; w, q)$ denote a probability measure on ${\mathcal
B}({\mathbb R}^3)$
indexed by two parameters, $w \in[0,1]$ and $q \in\{0,1\}$, where
$(w,q)$ denotes a realization of
$(W_{n}^\prime,Q_{n})$ for
\[
W_n^\prime:= \frac{\log V_{n-1}}{\log u} \quad\mbox{and}\quad
Q_n:= {\mathbf1}_{\{T_u < {n}\}}.
\]
Set $W_n = W_n^{\prime} {\mathbf1}_{\{W_n^{\prime} \in[0,1]\}} +(
W_n^{\prime} \wedge1){\mathbf1}_{\{W_n^{\prime} > 1\}}$. Note that $(W_n,
Q_n)$ is ${\mathfrak F}_{n-1}$ measurable. Let $\nu_n(\cdot) = \nu(\cdot; W_{n}, Q_{n})$ be a random measure derived from the measure
$\nu$. Observe that, conditioned on ${\mathfrak F}_{n-1}$, $\nu_n$ is a
probability measure. Now, we assume that the family of random measures
$\{ \nu_n (\cdot) \}\equiv\{\nu(\cdot; W_n, Q_n)\}$ satisfy the
following regularity condition:

\textit{Condition} (C$_0$): $\mu\ll\nu$ for each pair $(w,q) \in[0,1]
\times\{0,1\}$, and
\[
{\mathbf E}_{\mathfrak D} \biggl[ \log\biggl( \frac{d\mu}{d\nu}(Y_n;
W_{n},Q_{n}) \biggr) \bigg| W_{n} = w,
Q_n = q \biggr]
\]
is piecewise continuous as a function of $w$.

Let ${\mathfrak M}$ denote the class of measures $\{ \nu_n \}$ where
$\nu$ satisfies (C$_0$).
Thus, we consider a class of distributions where we shift \textit{all three}
members of the driving sequence $Y_n = (\log A_n,B_n,D_n)$ in some way,
allowing dependence on the history of the process
through the parameters $(w,q)$.

Now suppose that simulation is performed using a modification of our
main algorithm, where $Y_n \sim\nu_n$ for some collection $\bfnu:= \{
\nu_1,
\nu_2,\ldots\} \in{\mathfrak M}$. Let ${\mathcal E}_u^{(\bfnu)}$
denote the corresponding importance sampling estimator.
Let $\hat{\pi}_k$ denote an empirical estimate for $\pi$, as described
in the discussion of our main algorithm,
and let ${\mathcal E}_{u,1}^{(\bfnu)},\ldots, {\mathcal E}_{u,n}^{(\bfnu
)}$ denote simulated estimates for ${\mathcal E}_{u}^{(\bfnu)}$
obtained by repeating this algorithm,
but with $\{ \nu_n \}$ in place of the dual measure \eqref{DUAL-MEAS}.
Then it is easy to see, using the arguments of Theorem~\ref{MainTh-1}, that
%
%e2.16 #&#
\begin{equation}
\label{sr100} \lim_{k \to\infty} \lim_{n \to\infty} \hat{
\pi}_{k} ({\mathcal C}) \widehat{\mathcal E}_{u,n}^{(\bfnu)}
= {\mathbf P} \{ V > u \},
\end{equation}
where $\widehat{\mathcal E}_{u,n}^{(\bfnu)}$ denotes the average of $n$
simulated samples of ${\mathcal E}_{u}^{(\bfnu)}$ (and depends on~$k$); cf.~\eqref{Cons}. It remains to compare
the variance of these estimators, which is the subject of the next theorem.

%
%th2.4 #&#
\begin{Th}\label{THoptimality}
Assume that the conditions of Theorems~\ref{MainTh-1} and \ref{MainTh-2} hold.
Let $\nu$ be a probability measure on ${\mathcal B}({\mathbb R}^3)$
indexed by parameters $w \in[0,1]$ and $q \in\{0, 1\}$,
and assume that $\bfnu\in{\mathfrak M}$. Then for any initial state
$v \in{\mathcal C}$,
%
%e2.17 #&#
\begin{equation}
\label{opt-1} \liminf_{u \to\infty} \frac{1}{\log u} \log\bigl(
u^{2\xi} {\mathbf E}_{\nu} \bigl[ \bigl( {\mathcal
E}_{u}^{(\bfnu)} \bigr)^2 | V_0 = v
\bigr] \bigr) \ge0.
\end{equation}
Moreover, equality holds in~\eqref{opt-1} if and only if $\nu(\cdot; w,
0) = \mu_{\xi}$
and $\nu(\cdot; w, 1) = \mu$ for all $w \in[0,1]$. Thus, the dual
measure in~\eqref{DUAL-MEAS} is the \textit{unique} optimal simulation strategy
within the class ${\mathfrak M}$.
\end{Th}

%
% PROOFS OF CONSISTENCY AND EFFICIENCY
%

%s3 #&#
\section{Proofs of consistency and efficiency}\label{sec3}
We start with consistency.

%
% PROOF OF THEOREM 2.2
%

\begin{pf*}{Proof of Theorem~\ref{ConsistencyTh}}
Let $K_0:=0$, $K_n:= \inf\{ i > K_{n-1}\dvtx  V_i \in{\mathcal C} \}$,
\mbox{$n \in\pintegers$}, denote the successive return times of $\{ V_n \}$
to ${\mathcal C}$.
Set
\[
X_n = V_{K_n},\qquad n=0,1,\ldots.
\]
Then we claim that the stationary distribution of $\{ X_n \}$ is given
by $\gamma(E) = \pi(E)/\pi({\mathcal C})$, where $\pi$ is the
stationary distribution of $\{ V_n\}$.

Notice that $\{ X_n \}$ is $\varphi$-irreducible and geometrically
ergodic [cf. \citet{JCAV11}, Lemma 5.1].
Now set ${\mathfrak N}_n:= \sum_{i=1}^n {\mathbf1}_{\{ V_i \in{\mathcal
C} \}}$. Then by the law of large numbers for Markov chains,
%
%e3.1 #&#
\begin{equation}
\label{march22a} \qquad\pi(E) = \lim_{n \to\infty} \frac{{\mathfrak N}_n}{n}
\Biggl( \frac{1}{{\mathfrak N}_n} \sum_{i=1}^{{\mathfrak N}_n} {
\mathbf1}_{\{ X_i \in E\}} \Biggr) = \pi({\mathcal C})\gamma(E)\qquad\mbox{ a.s., }
E \in{\mathcal B}({\mathcal C}).
\end{equation}
Hence $\gamma(E) = \pi(E)/\pi({\mathcal C})$.

Next, we assert that
${\mathbf P} \{ V > u \} = \pi({\mathcal C}) {\mathbf E}_\gamma [ N_u ]$.
To establish this equality, again apply the law of large numbers for
Markov chains to obtain that
%
%e3.2 #&#
\begin{eqnarray}\label{pf215}
{\mathbf P} \{ V > u \}&:=& \pi\bigl( (u,\infty) \bigr)
\nonumber\\[-8pt]\\[-8pt]
&=& \lim
_{n\to\infty} \frac{1}{n} \Biggl\{ \sum
_{i=0}^{K_{{\mathfrak N}_n}-1} {\mathbf1}_{\{ V_i > u \}} + \sum
_{i= K_{{\mathfrak N}_n}}^n {\mathbf1}_{\{ V_i > u \}} \Biggr\}\qquad\mbox{a.s.}\nonumber
\end{eqnarray}

By the Markov renewal theorem [\citet{IIPNEN85}, Lem\-ma~6.2],
we claim that the last term on the right-hand side (RHS) of this
equation converges to zero a.s.
To see this, let $I(n)$ denote the last regeneration time occurring in
the interval $[0,n]$,
let $J(n)$ denote the first regeneration time occurring after time $n$,
let $\tau$ denote a typical regeneration time.
Then by Lemma~6.2 of \citet{IIPNEN85} and the geometric ergodicity of $\{V_n\}$,
\begin{equation}\label{pf217}
\lim_{n \to \infty} \mathbf{E} \bigl[ e^{\varepsilon(J(n)-I(n))} \bigr] =
\frac{1}{\mathbf{E}[\tau]}\mathbf{E} \bigl[ \tau e^{\varepsilon \tau} \bigr]
 < \infty,\qquad\mbox{some } \varepsilon > 0.
\end{equation}
Now by Nummelin's split-chain construction [\citet{EN84}, Section 4.4]
and by the definition of $K_{{\mathfrak N}_n}$,
$I(n) \le K_{{\mathfrak N}_n} \le n \le J(n)-1$.  Hence by a Borel--Cantelli argument,
\begin{equation}\label{pf216}
\frac{1}{n}\sum_{i=K_{{\mathfrak N}_n}}^n \mathbf{1}_{\{V_i > u\}} \to 0\qquad\mbox{a.s. as } n \to \infty.
\end{equation}

%Then by Lemma 6.2 of \citet{IIPNEN85},
%%
%%e3.3 #&#
%{\mathbf1}_{\{ V_i > u \}} = \frac{1}{{\mathbf E} [ \tau ]} {\mathbf E} \Biggl[
%_{i=0}^n {\mathbf1}_{\{ V_i > u \}} \Biggr) \Bigg|
%V_0 \sim\eta\Biggr].
%%
%Now by Nummelin's split-chain construction [\citet{EN84},
%Section~4.4] and by the definition of $K_{{\mathfrak N}_n}$,
%$I(n) \le K_{{\mathfrak N}_n} \le n \le J(n)-1$. Hence
%%
%%e3.4 #&#
%{\mathbf1}_{\{ V_i > u \}} \le\frac{1}{n} \sum_{i= I(n)}^{J(n)-1}
%{\mathbf1}_{\{ V_i > u \}} \searrow0 \qquad\mbox{as } n \to\infty.

Next consider the first term on the RHS of~\eqref{pf215}. Assume $V_0$
has distribution $\gamma$.
For any $n \in\pintegers$, set $N_{u,n} = \sum_{i=K_{n-1}}^{K_n-1}
{\mathbf1}_{\{ V_i > u \}}$
(namely, the number of exceedances
above level $u$ which occur over the successive cycles starting
from ${\mathcal C}$). Let
$S_n^N = N_{u,1} + \cdots+N_{u,n}$, $n \in\pintegers$. It can be seen
that $\{ (X_n, N_{u,n})\}$ is a positive
Harris chain and, hence,
by another application of the law of large numbers for Markov chains,
%
%e3.5 #&#
\begin{equation}
\label{pf214} {\mathbf E}_\gamma[N_u ] = \lim
_{n \to\infty} \frac{S^N_n}{n}:= \lim_{n \to\infty}
\frac{1}{n} \sum_{n=0}^{K_n-1} {\mathbf
1}_{\{ V_i > u \}} \qquad\mbox{a.s.}
\end{equation}
Since ${\mathfrak N}_n/n \to\pi({\mathcal C})$ as $n \to\infty$, it
follows from~\eqref{pf215},
\eqref{pf216} and~\eqref{pf214} that
%
%e3.6 #&#
\begin{equation}
\label{pf218} {\mathbf P} \{ V > u \} = \lim_{n\to\infty}
\frac{{\mathfrak N}_n}{n} \Biggl( \frac{1}{{\mathfrak N}_n} \sum
_{i=0}^{K_{{\mathfrak N}_n}-1}
{\mathbf1}_{\{ V_i > u \}} \Biggr) = \pi({\mathcal C}) {\mathbf E}_\gamma
[N_u ].
\end{equation}
Finally recall ${\mathcal E}_u:= N_u e^{-\xi S_{T_u}}{\mathbf1}_{\{T_u <K\}
}$ and hence by an elementary change-of-measure argument [as in (\ref{march22d})
below], we have ${\mathbf E}_{\gamma}[N_u]={\mathbf E}_{\mathfrak D}[{\mathcal
E}_u]$.

To complete the proof, it remains to show that
%
%e3.7 #&#
\begin{eqnarray}
\label{pf2111} \qquad\lim_{k \to\infty} {\mathbf E}_{\mathfrak D} \bigl[
N_u e^{-\xi S_{T_u}} {\mathbf1}_{\{ T_u < K \}} | V_0
\sim\hat{\gamma}_k \bigr]
&=& {\mathbf E}_{\mathfrak D} \bigl[ N_u
e^{-\xi S_{T_u}} {\mathbf1}_{\{ T_u < K \}} | V_0 \sim\gamma\bigr],
\end{eqnarray}
where $S_n:= \sum_{i=1}^n \log A_i$.
Set
%
%e3.8 #&#
\begin{equation}
\label{pf2112} H(v) = {\mathbf E}_{\mathfrak D} \bigl[ {\mathbf E}_{\mathfrak D}
[N_u | {\mathfrak F}_{T_u} ] e^{-\xi S_{T_u}} {\mathbf
1}_{\{ T_u < K \}} | V_0 = v \bigr].
\end{equation}
We now claim that $H(v)$ is uniformly bounded in $v \in{\mathcal C}$.
To establish this claim, first apply Proposition 4.1 of \citet{JCAV11}
to obtain that
%
%e3.9 #&#
\begin{equation}
\label{pf2113} {\mathbf E}_{\mathfrak D} [ N_u | {\mathfrak
F}_{T_u} ] {\mathbf1}_{\{T_u < K\}} \le\biggl( C_1(u)
\log\biggl(\frac{V_{T_u}}{u} \biggr) + C_2(u) \biggr) {\mathbf
1}_{\{ T_u < \tau\}},
\end{equation}
where $\tau\ge K$ is the first regeneration time and $C_i(u) \to C_i <
\infty$ as $u \to\infty$ ($i=1,2$).
Moreover, for $Z_n:= V_n/(A_1 \cdots A_n)$, we clearly have
%
%e3.10 #&#
\begin{equation}
\label{pf2114} e^{-\xi S_{T_u}} = u^{-\xi} \biggl( \frac{V_{T_u}}{u}
\biggr)^{-\xi} Z_{T_u}^\xi.
\end{equation}
Substituting the last two equations into~\eqref{pf2112} yields
%
%e3.11 #&#
\begin{equation}
\label{pf2110} \bigl|H(v)\bigr| \le\Theta{\mathbf E}_{\mathfrak D} \bigl[
\bigl|Z^\xi_{T_u} {\mathbf1}_{\{ T_u < \tau\}} \bigr| | V_0 =
v \bigr] \le\widebar{\Theta}
\end{equation}
for finite constants $\Theta$ and $\widebar{\Theta}$,
where the last step was obtained by \citet{JCAV11}, Lemma 5.5(ii).
Consequently, $H(v)$ is bounded uniformly in $v \in{\mathcal C}$.

Since
$\hat{\gamma}_k$ and $\gamma$ are both supported on ${\mathcal C}$, it
then follows since $\hat{\gamma}_k \Rightarrow\gamma$
that
\[
\lim_{k \to\infty} \int_{\mathcal C} H(v) \,d\hat{
\gamma}_k(v) = \int_{\mathcal C} H(v) \,d\gamma(v),
\]
which is~\eqref{pf2111}.
\end{pf*}

Before turning to the proof of efficiency, it will be helpful to have a
characterization of the return times of $\{ V_n \}$
to the set ${\mathcal C}$ when $Y_n \sim\mu_\beta$ for $\beta\in\operatorname{dom} (\Lambda)$, where $Y_n:= (\log A_n,B_n,D_n)$
and $\mu_\beta$ is defined according to~\eqref{meas-shift}.
First let
\[
\lambda_\beta(\alpha) = \int_{{\mathbb R}^3}
e^{\alpha x} \,d\mu_\beta(x,y,z),\qquad\Lambda_\beta(
\alpha) = \log\lambda_\beta(\alpha),\qquad\alpha\in{\mathbb R}
\]
and note by the definition of $\mu_\beta$ that
%
%e3.12 #&#
\begin{equation}
\label{MarkovLm-1} \Lambda_\beta(\alpha) = \Lambda(\alpha+ \beta) -
\Lambda(\beta).
\end{equation}

Recall that if $P$ denotes the transition kernel of $\{ V_n\}$, then we
say that $\{ V_n \}$ satisfies a \textit{drift condition} if there
exists a function $h\dvtx  {\mathbb R} \to[0,\infty)$ such that
{\renewcommand{\theequation}{${\mathcal D}$}
%e3.13 #&#
\begin{equation}\label{mmm}
\int_{\mathbb S} h(y) P(x,dy) \le\rho h(x)\qquad\mbox{for all } x
\notin{\mathcal C},
\end{equation}}\setcounter{equation}{12}%
where $\rho\in(0,1)$ and ${\mathcal C}$ is some Borel subset of
${\mathbb R}$.

%
% LEMMA 3.3
%

%
%le3.1 #&#
\begin{Lm}\label{le3.1}
Assume Letac's model~E, and suppose that \textup{(H$_1$)}, \textup{(H$_2$)} and \textup{(H$_3$)}
are satisfied.
Let $\{ V_n \}$ denote the forward
recursive sequence generated by this SFPE under the measure $\mu_\beta$,
chosen such that $\inf_{\alpha>0} \lambda_\beta(\alpha) < 1$.
Then the drift condition
(\ref{mmm}) holds with $h(x) = |x|^\alpha$, where $\alpha>0$ is
any constant satisfying the equation
$\Lambda_\beta(\alpha)<0$. Moreover, we may take $\rho= \rho_\beta$
and ${\mathcal C} = [-M_\beta,M_\beta]$,
where
%
%e3.13 #&#
\begin{equation}
\label{lm331n} \rho_\beta:= t \lambda_\beta(\alpha)\qquad\mbox{for
some } t \in\biggl(1,\frac{1}{\lambda_\beta(\alpha)} \biggr)
\end{equation}
and
%
%e3.14 #&#
\begin{equation}
\label{lm332n} \qquad M_\beta:= \cases{ \bigl( {\mathbf E}_\beta
\bigl[ \widetilde{B}^\alpha\bigr] \bigr)^{1/\alpha} \bigl(
\lambda_\beta(\alpha) (t-1 ) \bigr)^{-1/\alpha}, &\quad if $\alpha
\in(0,1)$,
\vspace*{5pt}\cr
\bigl( {\mathbf E}_\beta\bigl[ \widetilde{B}^{\alpha}
\bigr] \bigr)^{1/\alpha} \bigl( \bigl( \lambda_\beta(\alpha)
\bigr)^{1/\alpha} \bigl(t^{1/\alpha}-1 \bigr) \bigr)^{-1}, &
\quad if $\alpha\ge1$.}
\end{equation}
Furthermore, for any $(\rho_\beta,M_\beta)$ satisfying this pair of equations,
%
%e3.15 #&#
\begin{equation}
\label{lm333n} \sup_{v \in{\mathcal C}} {\mathbf P}_\beta\{ K > n |
V_0 = v \} \le\rho_\beta^{n}\qquad\mbox{for all } n \in\pintegers.
\end{equation}
\end{Lm}

\begin{pf}
Let $\widetilde{B}_n:= A_n |D_n| + |B_n|$.
If $\alpha\ge1$, then Minkowskii's inequality yields
%
%e3.16 #&#
%e3.17 #&#
\begin{eqnarray}\label{lmpf332}
&& {\mathbf E}_\beta \bigl[ |V_1|^\alpha
| V_0=v \bigr]\nonumber
\\
&&\qquad \le \bigl( \bigl({\mathbf E}_\beta
\bigl[A^\alpha\bigr] \bigr)^{1/\alpha} v + \bigl({\mathbf
E}_\beta\bigl[ \widetilde{B}^\alpha\bigr] \bigr)^{1/\alpha}
\bigr)^{\alpha}
\\
&&\qquad  = \rho_\beta v^\alpha\biggl( \frac{1}{t^{1/\alpha}} +
\frac{ ({\mathbf E}_\beta [ \widetilde{B}^\alpha ] )^{1/\alpha}}{\rho_\beta
^{1/\alpha}v} \biggr)^\alpha\qquad\mbox{where } \rho_\beta:= t \lambda _\beta(\alpha).\nonumber
\end{eqnarray}
Then (\ref{mmm}) is established. For $M_\beta$, set
$t^{-1/\alpha} + ( {\mathbf E}_\beta [ \widetilde{B}^\alpha ] )^{1/\alpha
}/(\rho_\beta^{1/\alpha}v) = 1$ and solve for $v$.
Similarly, if $\alpha< 1$, use $|x +y|^\alpha\le|x|^\alpha
+ |y|^\alpha$, $\alpha\in(0,1]$, in place of Minkowskii's inequality.
Then~\eqref{lm333n} follows by a standard argument, as in
\citet{EN84} or \citet{JCAV11}, Remark~6.2.
\end{pf}

We now introduce some additional notation which will be needed in the
proof of Theorem~\ref{MainTh-1}. Let $A_0 \equiv1$ and, for
any $n=0,1,2,\ldots,$ set
\begin{eqnarray*}
{\mathcal P}_n &=& A_0 \cdots A_n,\qquad S_n = \sum_{i=0}^n \log
A_i,
\\
Z_n &=& \frac{V_n}{A_0 \cdots A_n}\quad\mbox{and} \quad\widebar{Z}^{(p)} =
\sum_{n=0}^\infty
\frac{\widetilde{B}_n}{
A_0 \cdots A_n} {\mathbf1}_{\{ K > n \}},
\end{eqnarray*}
where
%
%e3.18 #&#
\begin{equation}
\label{defBtilde} \widetilde{B}_0 = |V_0| \quad\mbox{and} \quad\widetilde{B}_n = A_n |D_n| +
|B_n|.
\end{equation}

Also introduce the dual measure with respect to an arbitrary measure
$\mu_\alpha$, where $\alpha\in\operatorname{dom} (\Lambda)$.
Namely, define
{\renewcommand{\theequation}{${\mathfrak D}_\alpha$}
%e3.19 #&#
\begin{equation}
\label{defDUALgen} {\mathfrak L} (\log A_n,B_n,D_n
) = \cases{ \mu_{\alpha}, &\quad for $n=1,\ldots,T_u$,
\vspace*{2pt}\cr
\mu, &\quad for $n > T_u$.}
\end{equation}}\setcounter{equation}{17}%
Note that it follows easily from this definition that
for any r.v. ${\mathcal U}$ which is measurable with respect to ${\mathfrak F}_K$,
%
%e3.18 #&#
\begin{equation}
\label{march22d} {\mathbf E} [ {\mathcal U} {\mathbf1}_{\{T_u < K \}} ] = {\mathbf
E}_{\mathfrak D} \bigl[ \bigl(\lambda(\alpha) \bigr)^{T_u}
e^{-\alpha S_{T_u}} {\mathcal U} {\mathbf1}_{\{ T_u < K \}} \bigr],
\end{equation}
an identity which will be useful in the following.\vadjust{\goodbreak}

\begin{pf*}{Proof of Theorem~\ref{MainTh-1}}
Assume $V_0 =v \in{\mathcal C}$. We will show that the result holds
uniformly in
$v \in{\mathcal C}$.

\textit{Case} 1: $\lambda(\alpha)< \infty$, for some $\alpha< -\xi$.

To evaluate
\[
{\mathbf E}_{\mathfrak D} \bigl[ {\mathcal E}_u^2 \bigr]:= {\mathbf
E}_{\mathfrak D} \bigl[ N_u^2
e^{-2\xi S_{T_u}} {\mathbf1}_{\{ T_u < K \}} \bigr],
\]
first
note that
$V_n e^{-S_n}:= V_n/{\mathcal P}_n:= Z_n$.
Since $V_{T_u} > u$, it follows that $0 \le u e^{-S_{T_u}}
\le Z_{T_u}$.
Moreover, as in the proof of Lemma 5.5 of \citet{JCAV11} [cf.
(5.27), (5.28)], we obtain
\[
Z_n \le\sum_{i=0}^n
\frac{\widetilde{B}_i}{{\mathcal P}_i}\qquad\mbox{implying } Z_{T_u} {\mathbf
1}_{\{T_u < K \}}\le\sum_{n=0}^\infty
\frac{\widetilde{B}_n}{{\mathcal P}_n} {\mathbf1}_{\{ n \le T_u < K \}}.
\]
Consequently,
%
%e3.19 #&#
\begin{equation}
\label{pf221} u^{2\xi} {\mathbf E}_{\mathfrak D} \bigl[{\mathcal
E}_u^2 \bigr] \le{\mathbf E}_{\mathfrak D}
\Biggl[N_u^2 \Biggl( \sum_{n=0}^\infty
\frac{ \widetilde{B}_n}{{\mathcal P}_n} {\mathbf1}_{\{n \le T_u
< K \}} \Biggr)^{2\xi} \Biggr].
\end{equation}
If $2\xi\ge1$,
apply Minkowskii's inequality to the RHS to obtain
%
%e3.20 #&#
\begin{eqnarray}
\label{pf221a} \bigl( u^{2\xi} {\mathbf E}_{\mathfrak D} \bigl[ {\mathcal
E}_u^2 \bigr] \bigr)^{1/2\xi} & \le&\sum
_{n=0}^\infty\biggl( {\mathbf E}_{\mathfrak D} \biggl[
N_u^2 \biggl( \frac{\widetilde{B}_n}{{\mathcal P}_n} \biggr)^{2\xi}
{\mathbf1}_{\{ n\le T_u < K \}} \biggr] \biggr)^{1/2\xi}
\nonumber\\[-8pt]\\[-8pt]
& =& \sum_{n=0}^\infty\bigl( {\mathbf
E} \bigl[ N_u^2 {\mathcal P}_n^{-\xi}
\widetilde{B}_n^{2\xi} {\mathbf1}_{\{ n \le T_u < K \}} \bigr]
\bigr)^{1/2\xi},\nonumber
\end{eqnarray}
where the last step follows from~\eqref{march22d}.
Using the independence of $(A_n,\widetilde{B}_n)$ and ${\mathbf1}_{\{ n-1 <
T_u \wedge K \}}$, it follows
by an application of H\"{o}lder's inequality
that the left-hand side (LHS) of~\eqref{pf221a} is bounded above by
\[
\sum_{n=0}^\infty\bigl( {\mathbf E} \bigl[
N_u^{2r} \bigr] \bigr)^{1/2r\xi} \bigl( {\mathbf E}
\bigl[ \bigl( A_n^{-1} \widetilde{B}_n^2
\bigr)^{s\xi} \bigr] \bigr)^{1/2s\xi} \bigl({\mathbf E} \bigl[ {\mathcal
P}_{n-1}^{-s\xi} {\mathbf1}_{\{ n-1 < T_u \wedge K \}} \bigr]
\bigr)^{1/2s\xi},
\]
where $r^{-1} + s^{-1}=1$.
Set $\zeta= s \xi$ for the remainder of the proof.
The last term on the RHS of the previous equation may be expressed in
$\mu_{-\zeta}$-measure as
%
%e3.21 #&#
\begin{equation}
\label{march23aa} \qquad {\mathbf E} \bigl[ {\mathcal P}_{n-1}^{-\zeta} {\mathbf
1}_{\{ n-1 < T_u \wedge K \}} \bigr] = \bigl(\lambda(-\zeta) \bigr
)^{n-1} {\mathbf
P}_{-\zeta} \{ n-1 < T_u \wedge K \}.
\end{equation}
Substituting this last equation into the upper bound for~\eqref{pf221a},
we conclude that
%
%e3.22 #&#
\begin{equation}
\label{pf222}\qquad  \bigl( u^{2\xi} {\mathbf E}_{\mathfrak D} \bigl[ {\mathcal
E}_u^2 \bigr] \bigr)^{1/2\xi} \le\sum
_{n=0}^\infty{\mathcal J}_{n} \bigl( \bigl(
\lambda(-\zeta) \bigr)^{n-1} {\mathbf P}_{-\zeta} \{ n-1 <
T_u \wedge K \} \bigr)^{1/2\zeta},
\end{equation}
where
\[
{\mathcal J}_n:= \bigl( {\mathbf E} \bigl[ N_u^{2r}
\bigr] \bigr)^{1/2r\xi} \bigl( {\mathbf E} \bigl[ \bigl( A_n^{-1}
\widetilde{B}_n^2 \bigr)^{\zeta} \bigr]
\bigr)^{1/2\zeta},\qquad n=0,1,\ldots.
\]

Since $N_u \le K$, applying Lemma~\ref{le3.1} with $\beta=0$ yields
%
%e3.23 #&#
\begin{equation}
\label{june24} \sup_{v \in{\mathcal C}} {\mathbf E} \bigl[ N_u^{2r}
| V_0=v \bigr] < \infty\qquad\mbox{for any finite constant } r.
\end{equation}
Moreover, for sufficiently small $s > 1$ and $\zeta= s \xi$, it
follows by~\eqref{newcond-thm22} that
${\mathbf E} [ ( A^{-1} \widetilde{B}^2 )^{\zeta} ] < \infty$.
Thus, to show that the quantity on the LHS of~\eqref{pf222} is
finite, it suffices to show for some $\zeta> \xi$
and some $t > 1$,
%
%e3.24 #&#
\begin{equation}
\label{pf223} {\mathbf P}_{{-\zeta}} \{ n-1 < T_u \wedge K \}
\le\bigl( t \lambda(-\zeta) \bigr)^{-n+1}\qquad\mbox{for all } n \ge
N_0,
\end{equation}
where $N_0$ is a finite positive integer, uniformly in $u$ and
uniformly in $v \in{\mathcal C}$.

To this end, note that $\{ T_u \wedge K > n-1 \} \subset\{ K > n-1 \}$,
and by Lemma~\ref{le3.1} [using that $\min_\alpha\lambda_{-\zeta}(\alpha
) < (\lambda(-\zeta) )^{-1}$ by (\ref{MarkovLm-1})],
%
%e3.25 #&#
\begin{equation}
\label{june19} \sup_{v \in{\mathcal C}} {\mathbf P}_{-\zeta} \{ K > n-1 |
V_0= v \} \le\bigl(t \lambda(-\zeta) \bigr)^{-n+1},
\end{equation}
where ${\mathcal C}:= [-M, M]$ and $M > M_{-\xi}$. [Since $\zeta> \xi$
was arbitrary, we have replaced
$M_{-\zeta}$ with $M_{-\xi}$ in this last expression. We note that we
also require $M>M_0$ for~\eqref{june24} to hold.]
We have thus established
\eqref{pf223} for the case
$2\xi\ge1$.

If $2\xi< 1$, then the above argument can be repeated but using
the deterministic inequality $|x+y|^\alpha\le|x|^\alpha+ |y|^\alpha
$, $\alpha\in(0,1]$, in place of Minkowskii's inequality,
establishing the theorem for this case.

\textit{Case} 2: $\lambda(-\zeta)= \infty$ for $\zeta> \xi$, while ${\mathbf
E} [ (A^{-1}\widetilde{B})^\alpha ] < \infty$ for all $\alpha> 0$.

First assume $2\xi\ge1$.
Then, as before, $ ( u^{2 \xi}{\mathbf E}_{\mathfrak D} [ {\mathcal E}_u^2
] )^{1/2\xi}$ is bounded above by
the RHS of~\eqref{pf221a}. In view of the display following \eqref
{pf221a}, it is sufficient to show that uniformly in $v \in{\mathcal C}$
(for some set ${\mathcal C} = [-M,M])$,
%
%e3.26 #&#
\begin{equation}
\label{march23a} \sup_{n \in\pintegers} {\mathbf E} \bigl[{\mathcal
P}_{n-1}^{-\zeta}{\mathbf1}_{\{ n-1 < T_u \wedge K \}} \bigr] < \infty\qquad\mbox{for some } \zeta> \xi.
\end{equation}

Set $W_n = {\mathcal P}_{n-1}^{-\zeta}{\mathbf1}_{\{ n-1 < T_u \wedge K \}}$,
and first observe that ${\mathbf E} [W_n ] < \infty$. Indeed,
%
%e3.27 #&#
\begin{equation}
\label{march23b} |V_{n}| \le A_{n}|V_{n-1}| \biggl(
1 + \frac{\widetilde{B}_{n}}{A_{n}|V_{n-1}|} \biggr),\qquad n=1,2,\ldots
\end{equation}
and $n-1< T_u \wedge K \Longrightarrow|V_{i}| \in(M, u)$ for
$i=1,\ldots,n-1$.
Hence~\eqref{march23b} implies
%
%e3.28 #&#
%e3.29 #&#
\begin{eqnarray}\label{june19g}
A_i^{-\zeta} \le\biggl( \frac{u}{M}
\biggr)^\zeta\biggl(1 + \frac{\widetilde{B}_{i}}{M A_{i}} \biggr)^\zeta,
\nonumber\\[-8pt]\\[-8pt]
\eqntext{i=1,\ldots,n-1 \mbox{ on } \{ n-1 < T_u \wedge K \}.}
\end{eqnarray}
This equation yields an upper bound for ${\mathcal P}_{n-1}$.
Using the assumption that
${\mathbf E} [ (A^{-1} \widetilde{B})^\alpha] < \infty$ for all $\alpha> 0$,
we conclude by~\eqref{june19g} that
${\mathbf E} [W_n ] < \infty$.

Next let $\{ L_k \}$ be a sequence of positive real numbers such that
$L_k \downarrow0$ as $k \to\infty$, and set
$F_k = \bigcap_{i=1}^{k-1} \{A_i \ge L_k\}$.
Assume that $L_k$ has been chosen sufficiently
small such that
%
%e3.30 #&#
\begin{equation}
{\mathbf E} [W_k {\mathbf1}_{F_k^c} ] \le\frac{1}{k^2},\qquad k=1,2,\ldots.
\end{equation}
Then it suffices to show that
%
%e3.31 #&#
\begin{equation}
{\label{verify1}} \sum_{k=0}^\infty{\mathbf E} [
W_k {\mathbf1}_{F_k} ]<\infty.
\end{equation}

To verify (\ref{verify1}), set
$\widebar{A}_{0,k} = 1$ and introduce the truncation
\[
\widebar{A}_{n,k} = A_n {\mathbf1}_{\{A_n \ge L_k\}} +
L_k {\mathbf1}_{\{A_n < L_k\}},\qquad n=1,2,\ldots.
\]
Let $\lambda_k(\alpha)= {\mathbf E} [ \widebar{A}^\alpha_{1,k} ]$
and $\widebar{W}_k=(\widebar{A}_0 \cdots\widebar{A}_{k-1})^{-\zeta}{\mathbf1}_{\{ k-1 < T_u \wedge K \}}$. After a
change of measure [as in~\eqref{march22d},~\eqref{march23aa}], we obtain
%
%e3.32 #&#
\begin{equation}
\label{june19a} {\mathbf E} [ \widebar{W}_k ] \le\bigl(\lambda_{k}(-
\zeta)\bigr)^{k-1} {\mathbf E}_{-\zeta} [ {\mathbf1}_{\{ K > k-1 \}} {\mathbf
1}_{F_k} ].
\end{equation}

To evaluate the expectation on the RHS, start with the inequality
%
%e3.33 #&#
\begin{equation}
\label{june19forward} |V_{n,k}| \le\widebar{A}_{n,k}|V_{n-1,k}|
\biggl(1 + \frac{\widetilde{B}_{n}}{ \widebar{A}_{n,k}|V_{n-1,k}|} \biggr
),\qquad n=1,2,\ldots.
\end{equation}
Write ${\mathbf E}_{-\zeta,w} [ \cdot ] =
{\mathbf E}_{-\zeta} [\cdot|V_{0,k}= w ]$. Then for any $\beta> 0$, a
change of measure followed by an application of H\"{o}lder's
inequality yields
%
%e3.34 #&#
\begin{eqnarray}
\label{july14a} {\mathbf E}_{-\zeta,w} \bigl[ |V_{1, k}|^{\beta}
\bigr] & \le&\frac{w^\beta}{\lambda_k(-\zeta)} {\mathbf E} \biggl[ ( \widebar{A}_{1,k}
)^{\beta-\zeta} \biggl( 1 + \frac{\widetilde{B}_{1}}{w \widebar{A}_{1,k}} \biggr
)^\beta\biggr]
\nonumber\\[-8pt]\\[-8pt]
& \le& \rho_k w^\beta\biggl( t^{-q} {
\mathbf E} \biggl[ \biggl( 1 + \frac{\widetilde{B}_{1}}{ w\widebar{A}_{1,k}} \biggr
)^{q\beta} \biggr]
\biggr)^{1/q},\nonumber
\end{eqnarray}
where $\rho_k:= ( {\mathbf E} [ ( \widebar{A}_{1,k} )^{p(\beta-\zeta)} ]
)^{1/p} (t/\lambda_k(-\zeta) )$
and $p^{-1} + q^{-1} = 1$.

Set $\hat{\beta} = \argmin_\alpha\lambda(\alpha)$ and choose $\beta$
such that $p(\beta-\zeta) = \hat{\beta}$,
and assume that $p>1$ is sufficiently small such that $\rho_k< \infty$,
$\forall k$. Noting that $\lambda(\hat{\beta}) < 1$,
we conclude that for $t \in (1, ( \lambda(\hat{\beta}) )^{-1/p} )$ and
for some constant $\rho\in(0,1)$,
%
%e3.35 #&#
\begin{equation}
\label{july14bb} \qquad\lim_{k \to\infty} \lambda_k(-\zeta)
\rho_k:= t \lim_{k \to\infty} \bigl({\mathbf E} \bigl[ (
\widebar{A}_{1,k} )^{p(\beta-\zeta)} \bigr]\bigr)^{1/p} = t \bigl( \lambda(\hat{\beta})
\bigr)^{1/p} < \rho,
\end{equation}
where the second equality was obtained by observing that as $k \to
\infty$, $L_k \downarrow0$ and hence
$\lambda_k(\alpha) \downarrow\lambda(\alpha)$, $\alpha>0$.
Equation~\eqref{july14bb} yields that
$\lambda_k(-\zeta) \rho_k \le\rho$ for all \mbox{$k \ge k_0$}, and with
this value of $\rho$,~\eqref{july14a} yields
%
%e3.36 #&#
\begin{equation}
\label{june23d} {\mathbf E}_{-\zeta,w} \bigl[ |V_{1, k}|^{\beta}
\bigr] \le\frac{\rho w^\beta}{\lambda_k(-\zeta)}\qquad\mbox{for all }
k \ge k_0,
\end{equation}
provided that
%
%e3.37 #&#
\begin{equation}
\label{july13b} t^{-q}{\mathbf E} \biggl[ \biggl(1 + \frac{\widetilde{B}_{1}}{w\widebar{A}_{1,k}}
\biggr)^{q \beta} \biggr] \le1.
\end{equation}

Our next objective is to find a set ${\mathcal C} = [-M,M]$ such that
for all $w \notin{\mathcal C}$, \eqref{july13b}~holds. First assume $q
\beta\ge1$
and apply Minkowskii's inequality to the LHS of~\eqref{july13b}. Then
set this quantity equal to one, solve for $w$ and set $w = M_k$.
After some algebra, this yields
%
%e3.38 #&#
\begin{equation}
\label{june23ee} M_k = \frac{1}{t^{1/\beta}-1} \biggl({\mathbf E} \biggl[
\biggl( \frac{ \widetilde{B}_1}{\widebar{A}_{1,k}} \biggr)^{q\beta} \biggr]
\biggr)^{1/q \beta}.
\end{equation}
The\vspace*{1pt} quantity in parentheses tends to
${\mathbf E} [ ( A^{-1} \widetilde{B} )^{q\beta} ]$
as $k \to\infty$.
Using the assumption ${\mathbf E} [ (A^{-1} \widetilde{B} )^\alpha ] < \infty$
for $\alpha>0$, we conclude $M:= \sup_k M_k < \infty$.

If $q\beta<1$, then a similar expression is obtained for $M$ by using
the deterministic inequality $| x+y |^\beta\le|x|^\beta+ |y|^\beta$
in place of Minkowskii's inequality.

To complete the proof, iterate
\eqref{june23d} with ${\mathcal C} = [-M,M]$ (as in the proof of Lemma~\ref{le3.1}) to obtain that
%
%e3.39 #&#
\begin{equation}
\label{july13a} {\mathbf E}_{-\zeta} [ {\mathbf1}_{\{K > k-1\}} {\mathbf
1}_{F_k} ] \le\biggl( \frac{\rho}{\lambda_k(-\zeta)} \biggr
)^{-k+1}\qquad\mbox{for all } k \ge k_0.
\end{equation}
Note that on the set $F_k$, $\{ V_{n,k}\dvtx  1 \le n \le k \}$ and $\{ V_n\dvtx
1 \le n \le k \}$ agree, and thus $\{ K > k-1 \}$
coincides for these two sequences. Substituting~\eqref{july13a}
into~\eqref{june19a} yields~\eqref{verify1} as required. Finally, the
modifications needed when
$2\xi< 1$ follow along the lines of those outlined in case~1, so we
omit the details.
\end{pf*}

%The above theorems could have been proved under the weaker assumption
%that $\{ V_n \}$ is bounded above and below by
%processes having the form of Letac's Model E, namely, $V_n =
%f_n(V_{n-1})$, where
%A_n \max\{v,D_n^\ast\} + B_n^\ast\le f_n(v) \le A_n \max\{v,D_n\} +
%B_n,
%for an i.i.d. sequence $\{ (A_n,B_n,D_n, B_n^\ast, D_n^\ast)\dvtx  n=1,2,
%Indeed, by following the modifications outlined in
%by working entirely with the upper bound in~\eqref{newintro5}.
%Similarly, Theorem 2.1 can be established by noting that, under
%and possesses
%all of the required regularity properties we have used here. In this
%way, our algorithm
%can be extended to various other processes beyond Letac's Model E,
%such as the AR(1) process with ARCH(1) errors, as described in more
%detail in
%}

%
% EXAMPLES AND SIMULATIONS
%

%s4 #&#
\section{Examples and simulations}\label{sec4}
In this section we provide several examples illustrating the
implementation of our algorithm.

%s4.1 #&#
\subsection{The ruin problem with stochastic investments}\label{sec4.1}
Let the fluctuations in the insurance business be governed
by the classical Cram\'{e}r--Lundberg model,
%
%e4.1 #&#
\begin{equation}
\label{sr1} X_t = u + c t - \sum_{n=1}^{N_t}
\zeta_n,
\end{equation}
where $u$ denotes the company's initial capital, $c$ its premium income
rate, $\{ \zeta_n \}$ the claims losses,
and $N_t$ the number of Poisson claim arrivals occurring in $[0,t]$.
Let $\{ \zeta_n \}$ be i.i.d. and independent of $\{ N_t \}$. We now
depart from this classical model by assuming that
at discrete times $n=1,2,\ldots,$ the surplus capital is invested,
earning stochastic returns $\{ R_n \}$, assumed to be i.i.d. Let
$L_n:= - (X_n - X_{n-1} )$ denote the losses incurred by the insurance
business during the $n$th discrete time interval. Then
the total capital of the insurance company at time $n$ is described by
the recursive sequence of equations
%
%e4.2 #&#
\begin{equation}
\label{ex311} Y_n = R_n Y_{n-1} -
L_n,\qquad n=1,2,\ldots,\ Y_0 = u,
\end{equation}
where it is typically assumed that ${\mathbf E} [ \log R ] > 0$ and ${\mathbf
E} [ L ] < 0$.

Our objective is to estimate the probability of ruin,
%
%e4.3 #&#
\begin{equation}
\label{ex310} \psi(u):= {\mathbf P} \{ Y_n < 0, \mbox{ for some } n \in
\pintegers| Y_0 = u \}.
\end{equation}
By iterating~\eqref{ex311}, we obtain that
$Y_n = (R_1 R_2\cdots R_n) (Y_0 -{\mathcal L}_n )$,
where\break  \mbox{${\mathcal L}_n:= \sum_{i=1}^n L_i/(R_1 \cdots R_i)$}.
Thus $\psi(u) = {\mathbf P} \{{\mathcal L}_n > u$, some $n \}$.
Setting $ {\mathcal L} = ( \sup_{n \in{\mathbb Z}_+} {\mathcal L}_n )
\vee0$,
then by an elementary argument [as in \citet{JCAV11}, Section~3],
we obtain that
${\mathcal L}$ satisfies the SFPE
%
%e4.4 #&#
\begin{equation}
\label{cl10} {\mathcal L} \stackrel{\mathcal D} {=} ( A {\mathcal L} +
B )^+\qquad\mbox{where } A \stackrel{\mathcal D} {=} \frac{1}{R_1}
\mbox{ and } B \stackrel{\mathcal D} {=} \frac{L_1}{R_1}.
\end{equation}
This can be viewed as a special case of Letac's model~E with $D:=-B/A$.

Now take
%
%e4.5 #&#
\begin{equation}
\label{july10a} A_n = \exp\biggl\{- \biggl(\mu- \frac{\sigma^2}{2}
\biggr) - \sigma Z_n \biggr\}\qquad\mbox{for all }n,
\end{equation}
where $\{Z_n\}$ is an i.i.d. sequence of standard Gaussian r.v.'s. It
can be seen that
$\xi= 2\mu/\sigma^2 - 1$ and $\mu_\xi\sim\operatorname{Normal}(\mu-\sigma
^2/2,\sigma^2)$.\vspace*{2pt}

We set $\mu= 0.2$,
$\sigma^2 = 0.25$, $c=1$, $\{\zeta_n \} \sim\operatorname{Exp}(1)$ and
let $\{ N_t \}$ be a Poisson process
with parameter~$1/2$.

We implemented our algorithm to estimate the probabilities of ruin for
$ u=10, 100, 10^3, 10^4, 10^5$.
In all of our simulations, the distribution in step~1 was based on
$k=10^4$, and $V_{1000}$ was taken as an
approximation to the limit r.v. $V$. We arrived at this choice using
extensive exploratory analysis and two-sample comparisons using
Kolmogorov--Smirnov tests between $V_{1000}$ and other values of $V_n$,
where $n=2000$, 5000, 10,000 (with $p$-values $\ge0.185$). Also, it is
worthwhile to point out\vadjust{\goodbreak} here that by Sanov's theorem and Markov chain
theory, the difference between the approximating $V_{n^{\ast}}$ and $V$
on ${\mathcal C}$ is exponentially small, since ${\mathcal C}$ is in
the center of the distribution of $V$.

In implementing the algorithm, we chose $M=0$, since,
arguing as in the proof of Lemma~\ref{le3.1}, we obtain that $M_\beta= \min
_{i=1,2} M_\beta^{(i)}$, where
%
%e4.6 #&#
\begin{eqnarray}
\label{july11a} M_\beta^{(1)} &=& \inf_{\alpha\in(0,1)\cap\Phi}
\frac{\|B_1^+\|_{\beta,\alpha}}{(1-\|A_1\|_{\beta,\alpha}^\alpha
)^{1/\alpha}},
\nonumber\\[-8pt]\\[-8pt]
M_\beta^{(2)} &=& \inf _{\alpha\in[1,\infty)\cap\Phi} \frac{\|B_1^+\|_{\beta,\alpha
}}{1-\|A_1\|_{\beta,\alpha}}\nonumber
\end{eqnarray}
and $\Phi= \{ \alpha\in{\mathbb R}\dvtx  {\mathbf E}_\beta [ A^\alpha ] < 1
\}$.
(Here $\|\cdot\|_{\beta,\alpha}$ denotes the $L_{\alpha}$ norm
under the measure $\mu_\beta$.) As previously, we consider two cases,
$\beta= 0$ and $\beta= -\xi$.
For each of these cases, this infimum is computed numerically, yielding
$M_0 = 0 = M_{-\xi}$.

%t1 #&#
\begin{table}
\tabcolsep=0pt
\caption{Importance sampling estimation for the ruin probability with investments obtained\break  using $M=0$}\label{table1}
\begin{tabular*}{\tablewidth}{@{\extracolsep{\fill}}@{}lcccccc@{}}
\hline
$\bolds{u}$& ${\mathbf P} \bolds{\{V>u\}}$ & \textbf{LCL}& \textbf{UCL}& $\bolds{C}$ &\textbf{RE} & \textbf{Crude est.}\\
\hline
1.0e$+$01 & 5.86e$-$02 & 5.65e$-$02 & 6.07e$-$02 & 2.33e$-$01 & 1.84e$+$01&5.73e$-$02\\
% 2.0e$+$01 & 3.66e$-$02 & 3.52e$-$02 & 3.81e$-$02 & 2.21e$-$01 & 1.98e$+$01%&3.54e$-$02\\
1.0e$+$02 & 1.33e$-$02 & 1.28e$-$02 & 1.39e$-$02 & 2.11e$-$01 & 2.12e$+$01&1.29e$-$02\\
% 5.0e$+$02 & 4.95e$-$03 & 4.74e$-$03 & 5.15e$-$03 & 2.06e$-$01 & 2.09e$+$01%&4.85e$-$03\\
1.0e$+$03 & 3.27e$-$03 & 3.14e$-$03 & 3.41e$-$03 & 2.07e$-$01 & 2.12e$+$01&3.21e$-$03\\
% 5.0e$+$03 & 1.25e$-$03 & 1.19e$-$03 & 1.30e$-$03 & 2.06e$-$01 & 2.18e$+$01%&1.23e$-$03\\
1.0e$+$04 & 8.13e$-$04 & 7.78e$-$04 & 8.49e$-$04 & 2.04e$-$01 & 2.24e$+$01&8.01e$-$04\\
% 5.0e$+$04 & 3.06e$-$04 & 2.93e$-$04 & 3.20e$-$04 & 2.02e$-$01 & 2.22e$+$01%&3.27e$-$04\\
1.0e$+$05 & 1.98e$-$04 & 1.90e$-$04 & 2.07e$-$04 & 1.98e$-$01 & 2.16e$+$01&2.10e$-$04\\
\hline
\end{tabular*}
\end{table}

Table~\ref{table1} summarizes the probabilities of ruin (with $M=0$) and the
lower and upper
bounds of the 95\% confidence intervals (LCL, UCL) based on $10^6$
simulations. The confidence intervals in this and other examples
in this section are based on the simulations; that is, the lower 2.5\%
and upper 97.5\% quantiles of the simulated values of ${\mathbf P}\{V >u\}$.
We also evaluated the true constant $C(u):= {\mathbf P} \{ V > u \} u^\xi$
[which would appear in~\eqref{intro3} if this
expression were exact], and the relative error~(RE).
Even in the extreme tail---far below the probabilities of practical
interest in this problem---our algorithm works effectively and is
clearly seen to have bounded relative error. For comparison, we also
present the crude Monte Carlo estimates of the probabilities of ruin
based on $5 \times10^6$ realizations of $V_{2000}$.
We observe that for small values of $u$, the importance sampling
estimates and the crude Monte
Carlo estimates are close, \textit{which provides an empirical
validation of the algorithm for small values of $u$}.

%s4.2 #&#
\subsection{The $\mathrm{ARCH}(1)$ process}\label{sec4.2}
Now consider the ARCH(1) process, which models the squared returns on
an asset via the
recurrence equation
\[
R_n^2 = \bigl( a + b R_{n-1}^2
\bigr) \zeta_n^2 = A_n R_{n-1}^2
+ B_n,\qquad n=1,2,\ldots,
\]
where\vspace*{-1pt} $A_n= b\zeta_{n}^2$, $B_n = a \zeta_{n}^2$, and $\{ \zeta_n \}$
is an i.i.d. Gaussian sequence.
Setting \mbox{$V_n=R_n^2$}, we see that $V:= \lim_{n \to\infty} V_n$
satisfies the SFPE
$V \stackrel{\mathcal D}{=} A V + B$, and it is easy to verify that the
assumptions of our theorems are satisfied.
Then it is of interest to determine ${\mathbf P} \{ V > u \}$
for large $u$.

Next we implement our algorithm to estimate these tail probabilities.
As in the previous example, we identify $V_{1000}$ as an approximation
to $V$.
Turning to identification of $M$, recall that in the previous example,
we worked with a sharpened form of the formulas in Lemma~\ref{le3.1}; however,
in other examples, this approach may, like Lemma~\ref{le3.1}, yield a poor
choice for $M$. This is due to the fact that these types of estimate
for $V_n^{\alpha}$ typically use Minkowskii- or H\"{o}lder-type
inequalities, which are usually not very sharp. We now outline an
alternative method for obtaining $M$ and demonstrate that it yields
meaningful answers from a practical perspective. In the numerical
method, we work
directly with the conditional expectation and avoid upper-bound
inequalities. We emphasize that this procedure applies to any process
governed by Letac's model~E.

\subsubsection*{Numerical procedure for calculating $M$}
The procedure involves a Monte Carlo method for calculating the
conditional expectation appearing in the drift condition, that is,
for evaluating
\[
{\mathbf E}_\beta\biggl[ \biggl( \frac{V_1}{V_0} \biggr)^\alpha
\bigg| V_0=v \biggr] = {\mathbf E}_\beta\biggl[ \biggl( A \max
\biggl\{ \frac{D}{v},1 \biggr\} + \frac{B}{v} \biggr)^\alpha
\biggr],
\]
when $\beta=0$ and $\beta= -\xi$. The goal is to find an $\alpha$ such
that $M:= \max\{M_0,M_{-\xi}\}$ is minimized, where $M_\beta$
satisfies
\[
{\mathbf E}_\beta\biggl[ \biggl( A \max\biggl\{ \frac{D}{v},1
\biggr\} + \frac{B}{v} \biggr)^\alpha\biggr] \le
\rho_\beta\qquad\mbox{for all } v> M_\beta\mbox{ and some }\rho_\beta\in(0,1).
\]
In this expression, $\alpha$ is chosen such that ${\mathbf E}_\beta [
A^\alpha ] \in(0,1)$, and hence we expect that
$\rho_\beta\in ( {\mathbf E}_\beta [ A^\alpha ],1 )$. Note that $M_\beta
$ depends on the choice of $\alpha$; thus, we also minimize over
all possible $\alpha$ such that ${\mathbf E}_\beta [ A^\alpha ] \in(0,1)$.

Let $\{(A_i, B_i, D_i)\dvtx  1 \le i \le N\}$ denote a collection of i.i.d.
r.v.'s having the same distribution as $(A, B, D)$. Then the numerical
method for finding an optimal choice of $M$ proceeds as follows.

First, using a root finding algorithm such as Gauss--Hermite quadrature,
solve for $\xi$ in the equation ${\mathbf E}[A^\xi]=1$.
Next, for ${\mathbf E}_\beta [ A^\alpha ] < 1$, use a Monte Carlo procedure
with sample size $N$ to compute ${\mathbf E}_\beta [ |V_1|^{\alpha}| V_0= v ]$ and solve for $v$ in the formula
\[
\frac{1}{N} \sum_{i=1}^N \biggl
\llvert A_i \max\biggl\{ \frac{D_i}{v}, 1 \biggr\} +
\frac{B_i}{v}\biggr\rrvert^{\alpha} =\rho_\beta,
\]
where this quantity is computed in the $\beta$-shifted measure for
$\beta\in\{0,-\xi\}$ and where $\rho_\beta< 1$. Then select $\alpha$
so that it provides the smallest possible value of~$v$. Choose $M_\beta
> v$ for $\beta= 0$ and $\beta= -\xi$.
Finally, set $M= \max\{M_0, M_{-\xi}\}$.

\subsubsection*{Implementation}
We set $b=4/5$ and considered the values $a\dvtx  1.9\times10^{-5}, 1$.
It can be shown that
\[
{\mathbf E} \bigl[A_n^\alpha\bigr] = \frac{(2b)^\alpha\Gamma(\alpha
+1/2)}{\Gamma(1/2)}.
\]
We solved the equation ${\mathbf E} [A_n^\xi]=1$ using Gauss--Hermite
quadrature to obtain $\xi=1.3438$. Under the
$\xi$-shifted measure, $A_n=bX_n$ and $B_n=aX_n$, where $X_n\sim\Gamma
(\xi+1/2, 2)$.
Using the formulas in~\eqref{july11a} for $M$, we obtained [upon taking
the limit as $\delta\to0$ and using the Taylor approximation
$\Gamma(\delta+ 1/2) = \Gamma(1/2) + \delta\Gamma'(1/2) + O(\delta
^2)$] that
$M_0=0.362, 6.879\times10^{-6}$ when $a=1$, $1.9\times10^{-5}$, respectively.
Moreover, by applying the numerical method we have just outlined, it
can be seen that $M_{-\xi} = 0$.
[In contrast, by applying Lemma~\ref{le3.1} directly, one obtains $M_{-\xi} =
\infty$ since $\lambda(-\xi) = \infty$.]

%t2 #&#
\begin{table}
\tabcolsep=0pt
\caption{Importance sampling estimation for the tail probability of ARCH(1) financial process with
$a=1$, $1.9\times10^{-5}$}\label{table2}
\begin{tabular*}{\tablewidth}{@{\extracolsep{\fill}}@{}lcccccc@{}}\hline
$\bolds{u}$& $\bolds{\bP\{V > u\}}$ & \textbf{LCL}& \textbf{UCL}& $\bolds{C}$ &\textbf{RE}& \textbf{Crude est.}\\
\hline
\multicolumn{7}{@{}c@{}}{$a=1$}\\
1.0e$+$01 & 7.73e$-$02 & 7.64e$-$02 & 7.83e$-$02 & 1.71e$+$00 & 6.21e$+$00&7.75e$-$02\\
% 2.0e$+$01 & 3.43e$-$02 & 3.35e$-$02 & 3.51e$-$02 & 1.92e$+$00 & 1.18e$+$01%&3.43e$-$02\\
1.0e$+$02 & 4.34e$-$03 & 4.23e$-$03 & 4.45e$-$03 & 2.11e$+$00 & 1.29e$+$01&4.28e$-$03\\
% 5.0e$+$02 & 5.07e$-$04 & 4.96e$-$04 & 5.18e$-$04 & 2.15e$+$00 & 1.13e$+$01%&5.21e$-$04\\
1.0e$+$03 & 2.04e$-$04 & 1.99e$-$04 & 2.09e$-$04 & 2.20e$+$00 & 1.28e$+$01&2.07e$-$04\\
% 5.0e$+$03 & 2.32e$-$05 & 2.28e$-$05 & 2.36e$-$05 & 2.17e$+$00 & 8.08e$+$00%&2.00e$-$05\\
1.0e$+$04 & 9.00e$-$06 & 8.88e$-$06 & 9.12e$-$06 & 2.14e$+$00 & 6.83e$+$00&9.00e$-$06\\
%5.0e$+$04 & 1.07e$-$06 & 1.05e$-$06 & 1.10e$-$06 & 2.21e$+$00 & 1.27e$+$01%&2.00e$-$06\\
1.0e$+$05 & 4.11e$-$07 & 4.04e$-$07 & 4.18e$-$07 & 2.15e$+$00 & 8.51e$+$00 &NA
\\[3pt]
% 1.0e$+$01 & 1.62e$-$01 & 1.60e$-$01 & 1.64e$-$01 & 3.57e$+$00 & 5.99e$+$00
%&1.62e$-$01\\
% 2.0e$+$01 & 7.73e$-$02 & 7.64e$-$02 & 7.83e$-$02 & 4.33e$+$00 & 6.21e$+$00
%&7.78e$-$02\\
% 1.0e$+$02 & 1.08e$-$02 & 1.05e$-$02 & 1.11e$-$02 & 5.25e$+$00 & 1.34e$+$01
%&1.06e$-$02\\
% 5.0e$+$02 & 1.28e$-$03 & 1.25e$-$03 & 1.31e$-$03 & 5.43e$+$00 & 1.14e$+$01
%&1.33e$-$03\\
% 1.0e$+$03 & 5.07e$-$04 & 4.96e$-$04 & 5.18e$-$04 & 5.45e$+$00 & 1.13e$+$01
%&5.44e$-$04\\
% 5.0e$+$03 & 5.96e$-$05 & 5.81e$-$05 & 6.11e$-$05 & 5.57e$+$00 & 1.27e$+$01
%&7.70e$-$05\\
% 1.0e$+$04 & 2.32e$-$05 & 2.28e$-$05 & 2.36e$-$05 & 5.51e$+$00 & 8.08e$+$00
%&3.50e$-$05\\
% 5.0e$+$04 & 2.64e$-$06 & 2.60e$-$06 & 2.68e$-$06 & 5.44e$+$00 & 7.60e$+$00
%&3.00e$-$06\\
% 1.0e$+$05 & 1.07e$-$06 & 1.05e$-$06 & 1.10e$-$06 & 5.61e$+$00 & 1.27e$+$01
%&1.00e$-$06\\
\multicolumn{7}{@{}c@{}}{$a=1.9\times10^{-5}$}\\
1.0e$+$01 & 4.45e$-$08 & 4.38e$-$08 & 4.52e$-$08 & 9.82e$-$07 & 8.38e$+$00 &NA\\
% 2.0e$+$01 & 1.75e$-$08 & 1.72e$-$08 & 1.78e$-$08 & 9.80e$-$07 & 1.00e$+$01 &NA\\
1.0e$+$02 & 2.02e$-$09 & 1.98e$-$09 & 2.05e$-$09 & 9.82e$-$07 & 9.29e$+$00 &NA\\
%5.0e$+$02 & 2.66e$-$10 & 1.99e$-$10 & 3.32e$-$10 & 1.13e$-$06 & 1.27e$+$02 &NA\\
1.0e$+$03 & 9.59e$-$11 & 8.77e$-$11 & 1.04e$-$10 & 1.03e$-$06 & 4.38e$+$01 &NA\\
%5.0e$+$03 & 1.04e$-$11 & 1.02e$-$11 & 1.06e$-$11 & 9.75e$-$07 & 1.01e$+$01 &NA\\
1.0e$+$04 & 4.15e$-$12 & 4.05e$-$12 & 4.26e$-$12 & 9.85e$-$07 & 1.32e$+$01 &NA\\
%5.0e$+$04 & 4.78e$-$13 & 4.66e$-$13 & 4.91e$-$13 & 9.86e$-$07 & 1.34e$+$01 &NA\\
1.0e$+$05 & 1.91e$-$13 & 1.83e$-$13 & 1.99e$-$13 & 1.00e$-$06 & 2.19e$+$01 &NA\\
\hline
\end{tabular*}
\end{table}

%
%Table 4.2(a) Importance sampling estimation for the tail probability
%of ARCH(1) financial process with $a=1$.\\
%{\small
%$u$& $\bP(V > u)$ & LCL& UCL& $C$ &RE&Crude Est.\\
% 1.0e+01 & 7.73e-02 & 7.64e-02 & 7.83e-02 & 1.71e+00 & 6.21e+00
%&7.75e-02\\
% 2.0e+01 & 3.43e-02 & 3.35e-02 & 3.51e-02 & 1.92e+00 & 1.18e+01
%&3.43e-02\\
%1.0e+02 & 4.34e-03 & 4.23e-03 & 4.45e-03 & 2.11e+00 & 1.29e+01
%&4.28e-03\\
%5.0e+02 & 5.07e-04 & 4.96e-04 & 5.18e-04 & 2.15e+00 & 1.13e+01
%&5.21e-04\\
% 1.0e+03 & 2.04e-04 & 1.99e-04 & 2.09e-04 & 2.20e+00 & 1.28e+01
%&2.07e-04\\
%5.0e+03 & 2.32e-05 & 2.28e-05 & 2.36e-05 & 2.17e+00 & 8.08e+00
%&2.00e-05\\
%1.0e+04 & 9.00e-06 & 8.88e-06 & 9.12e-06 & 2.14e+00 & 6.83e+00
%&9.00e-06\\
%5.0e+04 & 1.07e-06 & 1.05e-06 & 1.10e-06 & 2.21e+00 & 1.27e+01
%&2.00e-06\\
% 1.0e+05 & 4.11e-07 & 4.04e-07 & 4.18e-07 & 2.15e+00 & 8.51e+00 &NA\\
%}
Table~\ref{table2} summarizes the simulation results for the tail probabilities
of the ARCH(1) process
based on $10^6$ simulations. We notice a substantial agreement between
the crude Monte Carlo estimates and those produced by our algorithm for
small values of $u$. More importantly, we observe that the relative
error remains bounded in all of the cases considered, while the
simulation results
using the state-dependent algorithm in
\citet{JBHLBZ11} \emph{show that the relative error based on
their algorithm increases
as the
parameter $u \to\infty$}.
When compared with the state-independent algorithm of \citet
{JBHLBZ11}, our simulations give comparable numerical results
to those they report, although direct comparison is difficult due to
the unquantified role of bias in their formulas.
(In contrast, from a numerical perspective, the bias is negligible in
our formulas, as it involves the convergence of a Markov chain near the
\textit{center} of its distribution, which is known to occur at a
geometric rate.)
We emphasize that our method also applies to a wider class of problems,
as illustrated by the previous example. Finally, we remark that a
variant of the ARCH(1) process is the GARCH$(1,1)$ financial process,
which can be implemented by similar methods. Numerical results for this
model are roughly analogous, but further complications arise which can
be addressed as in our preprint under the same title in Math arXiv. For
a further discussion of examples governed by Letac's model~E and its
generalizations, see \citet{JCAV11}, Section~3.

%s5 #&#
\section{Proofs of results concerning running time of the algorithm}\label{sec5}
The proof of the first estimate will rely on the following.

%
%le5.1 #&#
\begin{Lm}\label{le5.1}
Under the conditions of Theorem~\ref{MainTh-2}, there exist positive constants
$\beta$ and $\rho\in(0,1)$ such that
%
%e5.1 #&#
\begin{equation}
\label{june25} {\mathbf E}_{\xi} \bigl[h(V_n)|V_{n-1}
\bigr] \le\rho h(V_{n-1})\qquad\mbox{on } \{ V_{n-1} \ge
\widebar{M} \}
\end{equation}
for some $\widebar{M} < \infty$,
where $h(x):= x^{-\beta}{\mathbf1}_{ \{x > 1\}} + {\mathbf1}_{\{x \le1\}}$.
\end{Lm}

\begin{pf}
Assume without loss of generality (w.l.o.g.) that $V_{n-1}= v >1$. Then
by the strong Markov property,
\[
{\mathbf E}_{\xi} \bigl[h(V_n)|V_{n-1}=v \bigr]= {\mathbf
E}_{\xi} \bigl[V_1^{-\beta} {\mathbf
1}_{\{V_1 >1\}}|V_{0}= v \bigr] +{\mathbf P}_{\xi}
\{V_1 \le1 |V_{0}= v \}.
\]
Using assumption~\eqref{march22f}, we obtain that the second term on
the RHS is $o(v^{-\varepsilon})$,
while the first term can be expressed as
\[
v^{\beta}{\mathbf E}_{\xi} \bigl[\bigl(A_1 \max\bigl
\{v^{-1}D_1, 1\bigr\} + v^{-1}B_1
\bigr)^{-\beta}{\mathbf1}_{\{V_1 >1\}}|V_0= v \bigr] \sim
v^{\beta} {\mathbf E}_\xi\bigl[ A_1^{-\beta}
\bigr]
\]
as $v \to\infty$.
Next observe that ${\mathbf E}_{\xi} [A_1^{-\beta} ] = \lambda(\xi-\beta)<1$
if $0 < \beta< \xi$. Thus, choosing $ \beta= \varepsilon\in(0,\xi)$,
where $\varepsilon$ is given as in~\eqref{march22f}, we obtain
that the lemma holds for any $\rho= ({\mathbf E}_{\xi} [A_1^{-\varepsilon} ],1
) $ and $\widebar{M} < \infty$ sufficiently large.
\end{pf}

\begin{pf*}{Proof of Theorem~\ref{MainTh-2}}
We will prove~\eqref{thm23a}--\eqref{thm23c} in three steps, each
involving separate ideas and certain
preparatory lemmas.\vspace*{1pt}

\begin{pf*}{Proof of Theorem~\ref{MainTh-2}, step~1}
\textit{Equation}~\eqref{thm23a} \textit{holds.}
Let $\widebar{M}$ be given as in Lemma~\ref{le5.1}, and assume w.l.o.g. that $\widebar{M} \ge\max\{M,1\}$.
Let $L\equiv\sup\{n \in\pintegers\dvtx
V_n \in(-\infty, \widebar{M}] \}$ denote the last exit time of $\{ V_n \}
$ from $(-\infty,\widebar{M}]$.
Then it follows directly from the definitions that $K \le L$ on $\{ K<
\infty\}$, where we recall that $K$ is
the return time to the ${\mathcal C}$-set. Thus it is sufficient to
verify that ${\mathbf E}_\xi [ L ] < \infty$.

To this end, we introduce two sequences of random times. Set ${\mathcal
J}_0 = 0$
and \mbox{${\mathcal K}_0 = 0$} and, for each $i \in\pintegers$,
\[
{\mathcal K}_i= \inf\{n > {\mathcal J}_{i-1}\dvtx
V_n > \widebar{M} \}\quad\mbox{and} \quad{\mathcal
J}_i=\inf\bigl\{n> {\mathcal K}_i\dvtx  V_n
\in(-\infty, \widebar{M}] \bigr\}.
\]
Our main interest is in $\{ {\mathcal K}_i \}$, the successive times
that the process escapes from the interval $(-\infty,\widebar{M}]$, and
$\kappa_i:= {\mathcal K}_i - {\mathcal K}_{i-1}$.

Let ${\mathfrak N}$ denote the total number of times that $\{ V_n \}$
exits $(-\infty,\widebar{M}]$ and
subsequently returns to $(-\infty,\widebar{M}]$.
Then it follows that
\[
L < \sum_{i = 1}^{{\mathfrak N}+1} \kappa_i.
\]
Then by the transience of $\{ V_n \}$ in $\mu_\xi$-measure, it follows that
${\mathbf E}_\xi [ {\mathfrak N} ] < \infty$.\vadjust{\goodbreak}

It remains to show that ${\mathbf E}_\xi [ \kappa_i ] < \infty$,
uniformly in the starting state $V_{\kappa_{i-1}} \in(\widebar{M},\infty]$.
But note that the ${\mathbf E}_{\xi}[\kappa_i]$ can be divided into two
parts; first, the sojourn time that the process $\{ V_n \}$ spends in
$(\widebar{M},\infty)$
prior to returning to $(-\infty,\widebar{M}]$ and, second, the sojourn time
in the interval $(-\infty,\widebar{M}]$ prior to exiting again.
Now if $\widebar{K}$ denotes the first return time to $(-\infty,\widebar{M}]$,
then by Lemma~\ref{le5.1},
\[
{\mathbf P}_\xi\{ \widebar{K} = n | V_0 = v \} \le
\rho^n \frac{ h(v)}{h(\widebar{M})} \le\rho^n.
\]
Hence ${\mathbf E}_\xi [ \widebar{K} {\mathbf1}_{\{ \widebar{K} < \infty\}} | V_0 =
v ] \le\Theta< \infty$, uniformly in $v > \widebar{M}$.

Thus, to establish the lemma, it is sufficient to show that ${\mathbf
E}_\xi [ \widebar{N} | V_0 = v ] < \infty$,
uniformly in $v \in(-\infty,\widebar{M}]$,
where $\widebar{N}$ denotes the total number of visits of $\{ V_n \}$ to
$(-\infty, \widebar{M}]$. To this end, first note that
$[-\widebar{M},\widebar{M}]$ is petite. Moreover, it is easy to verify that
$(-\infty,-\widebar{M})$ is also petite for sufficiently
large~$\widebar{M}$. Indeed, for large~$\widebar{M}$ and $V_0<-\widebar {M}$, \eqref{intro1} implies
$V_1 = A_1 D_1 + B_1$ w.p. $p >0$. Thus, $\{ V_n \}$ satisfies a
minorization with small set $(-\infty,-\widebar{M})$.
Consequently $(-\infty, \widebar{M}]$ is petite and hence uniformly transient.
We conclude ${\mathbf E}_\xi [ \widebar{N} ] < \infty$, uniformly in $V_0 \in
(-\infty,\widebar{M}]$.
\end{pf*}

Before proceeding to step~2, we need a slight variant of Lemma~4.1 in
\citet{JCAV11}. In the following,
let $A^l$ be a typical ladder height of the process $S_n = \sum_{i=1}^n
\log A_i$
in its $\xi$-shifted measure.

%
%le5.2 #&#
\begin{Lm}\label{le5.2}
Assume the conditions of Theorem~\ref{MainTh-2}. Then
%
%e5.2 #&#
\begin{equation}
\label{prelimNR1} \lim_{u \to\infty} {\mathbf P}_\xi\biggl\{
\frac{V_{T_u}}{u} > y \bigg| T_u < K \biggr\} = {\mathbf P}_\xi
\{ \widehat{V} > y \}
\end{equation}
for some r.v. $\widehat{V}$, where for all $y \ge0$,
%
%e5.3 #&#
\begin{equation}
\label{prelimNR2} {\mathbf P}_\xi\{ \log\widehat{V} > y \} =
\frac{1}{{\mathbf E}_\xi [ A^l ]} \int_y^\infty{\mathbf
P}_\xi\bigl\{A^l > z \bigr\} \,dz.
\end{equation}
\end{Lm}

\begin{pf}
%Let $V_n^\prime= (A_1 \cdots A_n) X_n$, where $X_n = Z_n$ on $\{ Z_n
%> 0 \}$, and $X_n=1$ otherwise.
It can be shown that
%
%e5.4 #&#
\begin{equation}
\label{Lm41-new1} \frac{V_{T_u}}{u} \Rightarrow\widehat{V} \qquad\mbox{as }
u \to\infty
\end{equation}
in $\mu_\xi$-measure, independent of $V_0 \in{\mathcal C}$
[see \citet{JCAV11}, Lemma~4.1].

Set $y > 1$. Then by~\eqref{Lm41-new1},
${\mathbf P}_\xi \{ V_{T_u}/u > y \} \to{\mathbf P}_\xi \{ \widehat{V} > y \}$
as $u \to\infty$; and using the independence
of this result on its initial state, we likewise have that
$ {\mathbf P}_\xi \{ V_{T_u}/u > y | T_u \ge K \} \to{\mathbf P}_\xi \{ \widehat{V} > y \}$ as $u \to\infty$. Hence we conclude
\eqref{prelimNR1}, provided that $\liminf_{u \to\infty} {\mathbf P}_\xi
\{ T_u < K \} > 0 \}$.
But by the transience of $\{ V_n \}$, ${\mathbf P}_\xi \{ T_u < K \} \to
{\mathbf P}_\xi \{ K = \infty \}
>0$ as $u \to\infty$.
\end{pf}\noqed
\end{pf*}

\begin{pf*}{Proof of Theorem~\ref{MainTh-2}, step~2}
\textit{Equation}~\eqref{thm23b} \textit{holds.} With respect to the measure $\mu
_\xi$,
it follows by Lemma 9.13 of \citet{DS85} that
%
%e5.5 #&#
\begin{equation}
\label{lm414} \frac{T_u}{\log u} \to\frac{1}{\Lambda^\prime(\xi)}
\qquad\mbox{in probability}
\end{equation}
(since $\Lambda^\prime(\xi) = {\mathbf E}_\xi[\log A]$). Hence,
conditional on $\{ T_u < K \}$,
$ ( T_u/\log u ) \to (\Lambda^\prime(\xi) )^{-1}$ in probability.

To show that convergence in probability implies convergence in expectation,
it suffices to show that the sequence $\{ T_u/\log u \}$ is uniformly
integrable.
Let $\widebar{M}$ be given as in Lemma~\ref{le5.1}, and first suppose that $\widebar{M}
\le M$ and $\operatorname{supp} ( V_n) \subset[-M,\infty)$ for all $n$.
Then, conditional on $\{ T_u < K \}$,
\[
T_u > n\quad \Longrightarrow\quad V_i \in(\widebar{M},u),\qquad
i=1,\ldots,n.
\]
Now apply Lemma~\ref{le5.1}. Iterating~\eqref{june25}, we obtain ${\mathbf E} [
h(V_n) \prod_{i=1}^n {\mathbf1}_{V_i \notin{\mathcal C}} | V_0 ]
\le\rho^n h(V_0)$, $n=1,2,\ldots.$
Then, using the explicit form of the function $h$ in Lemma~\ref{le5.1}, we conclude
that with $\beta$ given as in Lemma~\ref{le5.1},
%
%e5.6 #&#
\begin{equation}
\label{june25a} {\mathbf P}_\xi\{ T_u > n |
T_u < K \} \le\biggl( \frac{1}{{\mathbf P}_\xi \{ T_u < K \}} \biggr)
\rho^n u^\beta\qquad\mbox{for all } n.
\end{equation}
Now ${\mathbf P}_\xi \{ T_u < K \} \downarrow\Theta> 0$ as $ u \to\infty
$. Hence, letting
${\mathbf E}_\xi^{(u)} [ \cdot ]$ denote the expectation conditional on $
\{T_u < K \}$, we obtain that for some $\widebar{\Theta} < \infty$,
%
%e5.7 #&#
\begin{equation}
\label{lm4125} {\mathbf E}_\xi^{(u)} \biggl[ \frac{T_u}{\log u};
\frac{T_u}{\log u} \ge\eta\biggr] \le\widebar{\Theta}\rho^{\eta\log u}
u^\beta
\end{equation}
and for sufficiently large $\eta$, the RHS converges to zero as $u \to
\infty$.
Hence $\{ T_u/\log u \}$ is uniformly integrable.

If the assumptions at the beginning of the previous paragraph are not
satisfied, then write $T_u = L + (T_u - L)$, where
$L$ is the last exit time from the interval $(-\infty,\widebar{M}]$, as
defined in the proof of Theorem~\ref{MainTh-2}, step~1. Then $(T_u - L)$ describes
the length of the last excursion to level $u$ after exiting $(-\infty,\widebar{M}]$ forever. By a repetition of the argument just given, we
obtain that
\eqref{june25a} holds with $(T_u-L)$ in place of $T_u$; hence $\{ (T_u
- L)/\log u \}$ is uniformly integrable.
Next observe by the proof of Theorem~\ref{MainTh-2}, step~1, that ${\mathbf E}_\xi [
L/\log u ] \downarrow0$ as $u \to\infty$.
The result follows.
\end{pf*}

%
% LEMMA 5.2
%

%As a consequence of the techniques of the previous result, we also
%obtain:

%Under the conditions of the previous proposition, we have that
%${\mathbf E}_\xi\left[ K {\mathbf1}_{\{ K < \infty\}} \right] < \infty$.

%First assume $V_0 \ge M$, where ${\mathcal C}=[-M,M]$.
% Then repeat the previous argument with $u = M$ to obtain, similarly
%to~\eqref{lm4123},
%{\mathbf P}_\xi\left( \frac{S_n}{n} \le\log t \right) \le C_1 e^{-\theta
%n},\qquad\mbox{\rm for all} n,
%for certain positive constants $C_1$ and $\theta$. Combining
%this estimate with~\eqref{lm4122} (with $u=M$) yields
%{\mathbf P}_\xi\left( V_n \le1 \right) \le C_2 e^{-\theta^\prime n}
%for certain positive constants $C_2$ and $\theta^\prime$. The required
%result follows, since
%${\mathbf P} \left( K {\mathbf1}_{\{ K < \infty\}} \right) \le\sum_n n {\mathbf
%P} \left( V_n \le M \right) < \infty$.
%Finally, if $V_0 < M$ then we can define ${\mathcal K} = \inf\{n: V_n
%$K = {\mathcal K}+(K-{\mathcal K})$, and ${\mathbf E} \left[{\mathcal K}

%
% PROPOSITION 5.2
%

Turning now to the proof of the last equation in Theorem~\ref{MainTh-2}, assume
for the moment that
$(V_{0}/u) =v > 1$ (we will later remove this assumption); thus, the
process starts above level $u$ and so its dual measure agrees
with its initial measure. Also define
\[
L(z) = \inf\bigl\{ n\dvtx  |V_n| \le z \bigr\}\qquad\mbox{for any } z \ge0.
\]

%
% LEMMA 5.3
%

%
%le5.3 #&#
\begin{Lm}\label{le5.3}
Let
$(V_{0}/u) =v > 1$ and $t \in(0,1)$. Then under the conditions of
Theorem~\ref{MainTh-2},
%
%e5.8 #&#
\begin{equation}
\label{lm422} \lim_{u \to\infty} \frac{1}{\log u} {\mathbf E} \biggl[
L \bigl( u^t \bigr) \bigg| \frac{V_0}{u} = v \biggr] =
\frac{1-t}{|\Lambda^\prime(0)|}.
\end{equation}
\end{Lm}

\begin{pf}
For notational simplicity, we will suppress the conditioning on
$(V_0/u)=v$ in the proof.
We begin by establishing an upper bound. Define
\[
S_n^{(u)}:= \sum_{i=1}^n
X_i^{(u)} \qquad\mbox{where }X_i^{(u)}:= \log\bigl(A_i +
u^{-t} \bigl( A_i |D_i| + |B_i| \bigr)
\bigr).
\]
Then it can be easily seen that
%
%e5.9 #&#
\begin{equation}
\label{pf425} \log|V_n| - \log(vu) \le S_n^{(u)}\qquad\mbox{for all } n < L\bigl(u^t\bigr).
\end{equation}
Now let $\widetilde{L}_u(u^t) = \inf \{ n\dvtx  S_n^{(u)} \le-(1-t) \log u -
\log v \}$. Then
$L(u^t) \le\widetilde{L}_u(u^t)$ for all $u$.

By Wald's identity,
${\mathbf E} [S_{\widetilde{L}_u(u^t)} ] = {\mathbf E} [ X_1^{(u)} ] {\mathbf E} [
\widetilde{L}_u(u^t) ]$. Thus, letting
\[
O_u:= \bigl|S_{\widetilde{L}_u(u^t)} - (1-t) \log u - \log v\bigr|
\]
denote the overjump of $ \{S_n^{(u)} \}$ over a boundary
at level $(1-t)\log u + \log v$, we obtain
%
%e5.10 #&#
\begin{equation}
\label{pf426} L\bigl(u^t\bigr) \le\frac{(1-t) \log u + \log v + {\mathbf
E} [ O_u ]}{ |{\mathbf E} [ X_1^{(u)} ] |}.
\end{equation}
Since ${\mathbf E} [ X_1^{(u)} ] \to\Lambda^\prime(0)$ as $u \to\infty$,
the required upper bound will be established
once we show that
%
%e5.11 #&#
\begin{equation}
\label{pf426a} \lim_{u \to\infty} \frac{1}{\log u} {\mathbf E} [
O_u ] = 0.
\end{equation}

To establish~\eqref{pf426a}, note as in the proof of Lorden's
inequality [\citet{SA03}, Proposition V.6.1] that
${\mathbf E} [ O_u ] \le{\mathbf E} [ Y^2_u ]/ {\mathbf E} [ Y_u ]$, where $Y_u$
has the negative ladder height
distribution of the process $\{ S_n^{(u)} \}$. Next observe by
Corollary VIII.4.4 of \citet{SA03} that
%
%e5.12 #&#
\begin{equation}
\label{eq512} {\mathbf E} [ Y_u ] = m_u^{(1)}
e^{{\mathfrak S}_u} \to{\mathbf E} [ Y ] \qquad\mbox{as } u \to\infty,
\end{equation}
where $Y$ has the negative ladder height distribution of $\{ S_n \}$, and
$m_u^{(j)}:= | {\mathbf E} [ X^{(u)} ] |$, $j = 1,2,\ldots$ and
${\mathfrak S}_u:= \sum_{n=1}^\infty n^{-1} {\mathbf P} \{ S_n^{(u)} > 0 \}$.
We observe that ${\mathfrak S}_u$ is the so-called Spitzer series.
Similarly, an easy calculation [cf. \citet{DS85}, page~176] yields
%
%e5.13 #&#
\begin{equation}
\qquad\quad {\mathbf E} \bigl[ Y_u^2 \bigr] = m_u^{(2)}
e^{{\mathfrak S}_u} - 2 m_u^{(1)} e^{{\mathfrak S}_u} \sum
_{n=1}^\infty\frac{1}{n} {\mathbf E} \bigl[
\bigl(S_n^{(u)} \bigr)^+ \bigr] \to{\mathbf E} \bigl[
Y^2 \bigr],\qquad u \to\infty.
\end{equation}
Since ${\mathbf E} [ ( \log A )^3 ] < \infty\Longrightarrow{\mathbf E} [ Y^j
] < \infty$ for $j=1,2$, it follows that
${\mathbf E} [ O_u ] \to{\mathbf E} [Y^2 ]/{\mathbf E} [Y ] < \infty$, implying
\eqref{pf426a}. Thus~\eqref{lm422} holds as an upper bound.

To establish a corresponding lower bound, fix $s \in(t,1)$ and define
\[
\widetilde{L}\bigl(u^s\bigr) = \inf\bigl\{ n\dvtx  S_n \le
-(1-s) \log u - \log v \bigr\}.
\]
Observe that $V_n \ge A_n V_{n-1} - |B_n|$ for
$V_{n-1} \ge0$, and iterating yields
%
%e5.14 #&#
\begin{equation}
\label{prop4210} \qquad V_n \ge(A_1 \cdots A_n )
V_0 - W\qquad\mbox{where } W:= \lim_{n \to\infty}
\sum_{i=1}^{n} \prod
_{j=i+1}^n A_j |B_i|.
\end{equation}
Since $(V_0/u) = v$, it follows from the definition of $\widetilde{L}$ that
\[
\widetilde{L}\bigl(u^s\bigr) \ge n\quad \Longleftrightarrow\quad
(A_1 \cdots A_k)V_0 > u^s\qquad\mbox{for all } k < n.
\]
But by~\eqref{prop4210}, $(A_1 \cdots A_k)V_0 > u^s \Longrightarrow V_k
> u^t$ on $\{ W \le(u^s - u^t) \}$. Thus for all~$n$,
$\widetilde{L}(u^s) \ge n \Longrightarrow L(u^t) \ge n$ on $\{ W \le(u^s -
u^t) \}$, and consequently
%
%e5.15 #&#
\begin{equation}
\label{prop4211a} {\mathbf E} \bigl[ L\bigl(u^t\bigr) \bigr] \ge{\mathbf E}
\bigl[ \widetilde{L}\bigl(u^s\bigr); W \le\bigl(u^s -
u^t\bigr) \bigr].
\end{equation}

Next recall that for some  $\widebar{C}>0$,
%
%e5.16 #&#
\begin{equation}
\label{prop4211} {\mathbf P} \bigl\{ W > u^s - u^t \bigr\}
\sim\widebar{C} u^{-s\xi} \qquad\mbox{as } u \to\infty.
\end{equation}
As $\widetilde{L}(u^t)$ is the time required for the negative-drift random
walk $\{ S_n + \log v \}$
to reach the level $-(1-s)\log u$, Heyde's (\citeyear{CH66}) a.s. convergence
theorem for renewal processes gives that
%
%e5.17 #&#
\begin{equation}
\label{prop4212} \frac{\widetilde{L}(u^s)}{\log u} \to\frac{(1-s)}{|\Lambda
^\prime(0)|} \qquad\mbox{a.s. as } u \to\infty
\end{equation}
(since ${\mathbf E} [
\log A ] = \Lambda^\prime(0) < 0$). Hence
for any $\varepsilon> 0$,
%
%e5.18 #&#
\begin{equation}
\label{prop4213} \lim_{u \to\infty} {\mathbf P} \biggl\{ \frac{\widetilde{L}(u^t)}{\log u}
\notin(r-\varepsilon,r+\varepsilon) \biggr\} = 0\qquad\mbox{where } r:=
\frac{1-s}{|\Lambda^\prime(0)|}.
\end{equation}
Substituting~\eqref{prop4211} and~\eqref{prop4213} into \eqref
{prop4211a} and letting $\varepsilon\to0$, we obtain
%
%e5.19 #&#
\begin{equation}
\label{pf4214} \liminf_{u \to\infty} \frac{1}{\log u} {\mathbf E}
\bigl[ L \bigl( u^t \bigr) \bigr] \ge\frac{1-s}{|\Lambda^\prime(0)|}.
\end{equation}
The required lower bound follows by letting $s \downarrow t$.
\end{pf}
%
%BEGIN NEW CORRECTIONS FROM JEFF RECEIVED ON 05052013
%
% LEMMA 5.4
%

%
%le5.4 #&#
\begin{Lm}\label{le5.4}
Assume the conditions of the previous lemma. Then
%
%e5.20 #&#
\begin{equation}
\label{lm4220} \lim_{t \downarrow0} \biggl\{ \limsup
_{u \to\infty} \frac{1}{\log u} {\mathbf E} \bigl[ L(M) - L \bigl(
u^t \bigr) \bigr] \biggr\} = 0.
\end{equation}
\end{Lm}

\begin{pf}
Apply Lemma~~\ref{le3.1} with $\beta=0$ to obtain that, for some $\alpha> 0$,
\[
{\mathbf E} \bigl[ |V_n|^\alpha| V_{n-1} = w \bigr]
\le\rho|w|^\alpha\qquad\mbox{for all } w \notin{\mathcal C},
\]
where $\rho\in(0,1)$ and ${\mathcal C} = [-M,M]$, for some positive
constant $M$. Since this equation holds for all $n < L(M)$ (the first entrance
time into the set ${\mathcal C}$), iterating this equation yields
%
%e5.21 #&#
\begin{equation}
\label{lm4221} {\mathbf E} [ {\mathbf1}_{\{ L(M) > n \}} | V_0=w ] \le
\rho^n \biggl( \frac{|w|}{M} \biggr)^\alpha\qquad\mbox{for all } n.
\end{equation}

Now apply this equation to obtain an estimate for $L(M) - L(u^t)$. Since
$\llvert V_{L(u^t)}\rrvert \le u^t$, the previous equation [with
$V_{L(u^t)}$ in place of $V_0$] gives
%
%e5.22 #&#
\begin{equation}
\label{lm4222} {\mathbf P} \bigl\{ L(M) - L\bigl(u^t\bigr) > n \bigr\}
\le\rho^n \biggl( \frac{u^{t}}{M} \biggr)^\alpha\qquad\mbox{for all } n.
\end{equation}
Set $J_t(u) = L(M)- L(u^t)$ and $t^\prime= t\alpha/(-\log\rho) $.
Summing~\eqref{lm4222} over all $n \ge t^\prime\log u$ yields that
%
%e5.23 #&#
\begin{equation}
\label{lm4223} {\mathbf E} \bigl[ J_t(u) {\mathbf1}_{\{J_t(u) \ge t^\prime
\log u\}} \bigr]
\le\frac{\rho^{t^\prime\log u}}{1-\rho} \biggl( \frac{u^{t}}{M} \biggr
)^\alpha=
\frac{1}{(1-\rho)M^\alpha}.
\end{equation}
Hence
%
%e5.24 #&#
\begin{equation}
\label{lm4224} \limsup_{u \to\infty} \frac{1}{\log u} {\mathbf E}
\bigl[ L(M) - L \bigl( u^t \bigr) \bigr] \le t^\prime.
\end{equation}
Since $t^\prime\downarrow0$ as $t \downarrow0$,
we conclude~\eqref{lm4220}.
\end{pf}

%
% PROOF OF PROPOSITION 5.2
%

\begin{pf*}{Proof of Theorem~\ref{MainTh-2}, step~3}
\textit{Equation}~\eqref{thm23c} \textit{holds.}
By Lemmas \ref{le5.3}~and~\ref{le5.4},
%
%e5.25 #&#
\begin{equation}
\label{lm4230} H_u(v):= \frac{1}{\log u} {\mathbf E} \biggl[ L(M) \bigg|
\frac{V_0}{u} = v \biggr] \to\frac{1}{|\Lambda^\prime(0)|}\qquad\mbox
{as } u \to
\infty.
\end{equation}

Let $\hat{\mu}_u$, $\hat{\mu}$ denote the probability laws of the
r.v.'s $V_{T_u}/u$, $\widehat{V}$ appearing in the statement of Lemma~\ref{le5.2}.
Then, using the strong Markov property, it follows that $L(M)$,
conditional on $V_0/u \sim\hat{\mu}_u$, is equal in distribution to
$K - T_u$, conditional on $\{ T_u < K \}$. Thus it is sufficient to
verify that
%
%e5.26 #&#
\begin{equation}
\label{lm4231} \quad\qquad\lim_{u \to\infty} \frac{1}{\log u} {\mathbf E} \biggl[
L(M) \bigg| \frac{V_0}{u} \sim\hat{\mu}_u \biggr]:= \lim
_{u \to\infty}\int_{v \ge0} H_u(v) \,d
\hat{\mu}_u(v) = \frac{1}{|\Lambda^\prime(0)|}.
\end{equation}
This result will follow from~\eqref{lm4230}, provided that we can show
that the limit can be taken inside the integral in the above equation.

To do so, express the inner quantity in~\eqref{lm4231} as
%
%e5.27 #&#
\begin{equation}
\label{july1a} \int_{v \ge0} H_u(v)\,d ( \hat{
\mu}_u - \hat{\mu} ) (v) + \int_{v \ge0}
H_u(v) \,d\hat{\mu}(v).
\end{equation}
To deal with the first term, begin by obtaining an upper bound for $H_u(v)$.
First note by a slight modification of~\eqref{lm4223} [with $t=1$,
$\theta= -\log\rho$ and $J_1(u)$, $t^\prime$ replaced with $L(M)$,
$r$, resp.]
that
%
%e5.28 #&#
\begin{equation}
\label{lm4232} {\mathbf E} \bigl[ L(M) {\mathbf1}_{\{L(M) > r\log u\}} | |V_0|
\le u \bigr] \le\frac{u^{-r\theta}}{1-\rho} \biggl( \frac{u}{M}
\biggr)^\alpha
\end{equation}
for all $r >0$ and some $\alpha> 0$.
Now choose $r >\alpha/\theta$. Then the RHS is
bounded above by $\Theta_1 < \infty$, independent of $u$.
Consequently,
%
%e5.29 #&#
\begin{equation}
\label{pf42a1} \frac{1}{\log u} {\mathbf E} \bigl[ L(M) | |V_0| \le
u \bigr] \le r + \frac{\Theta_1}{\log u}.
\end{equation}

Next, we extend this estimate to the case where $(V_0/u) = v > 1$.
To this end, viewing an excursion time as the sum of the time to first
reach $[-u,u]$ and then reach ${\mathcal C}$,
we obtain
%
%e5.30 #&#
\begin{equation}
\label{pf42a2} \qquad {\mathbf E} \biggl[ L(M) \bigg| \frac{V_0}{u} = v \biggr] \le\sup
_{w \in(M,u]} {\mathbf E} \bigl[ L(M) | |V_0| =w \bigr] + {
\mathbf E} \biggl[ L(u) \bigg| \frac{V_0}{u} = v \biggr].
\end{equation}
For the second term, observe
\[
|V_{n-1}| > u\quad \Longrightarrow\quad |V_n| \le|V_{n-1}|
\biggl(A_n + \frac{\widetilde{B}_n}{u} \biggr);
\]
thus, $ {\mathbf E} [ L(u) | (V_0/u) = v ]$ is bounded above by the length
of time for the classical
random walk
\[
S_n^{(u)}:= S_{n-1}^{(u)} + \log\biggl(
A_n + \frac{\widetilde{B}_n}{u} \biggr),\qquad n=1,2,\ldots,
\]
starting from $S_0^{(u)} = \log(vu)$,
to reach the level $\log u$. Denote this sojourn time by $L^\ast(u)$.
Applying Lorden's inequality [\citet{SA03}, Proposition V.6.1] to
$ \{ S_n^{(u)} \}$,
we obtain [with $\Lambda^{\prime\prime\prime}(0) < \infty$] that
\[
{\mathbf E} \bigl[ L^\ast(u) \bigr] \le\Theta_2(u) \log v +
\Theta_3(u) \to\frac{\log v}{m_1}+ \frac{m_2}{m_1^2},\qquad u \to
\infty,
\]
where $m_i$ denotes the $i$th moment
of the ladder height distribution for the sequence $\{ \log A_i \}$;
cf. the discussion
following~\eqref{eq512} above.
Substituting this last bound and~\eqref{pf42a1} into~\eqref{pf42a2}, we
deduce that for some constant $\widebar{\Theta}$, uniformly in $u \ge
u_0$ for some finite constant $u_0$,
%
%e5.31 #&#
\begin{equation}
\label{pf42a3} H_u(v):= \frac{1}{\log u} {\mathbf E} \biggl[ L(M) \bigg|
\frac{V_0}{u} = v \biggr] \le\widebar{\Theta} + \frac{2\log v}{m_1}.
\end{equation}

Returning to~\eqref{july1a} and
using the above upper bound, we now show that
%
%e5.32 #&#
\begin{equation}
\label{july1b} %\left| \int_{v \ge0} H_u(v)\,d\left( \hat{\mu}_u - \hat{
\biggl\llvert\int_{v \ge0} \biggl(
\widebar{\Theta} + \frac{2 \log v}{m_1} \biggr) \,d ( \hat{\mu}_u - \hat{
\mu} ) (v) \biggr\rrvert\to0\qquad\mbox{as } u \to\infty.
\end{equation}
Since $\hat{\mu}_u \Rightarrow\hat{\mu}$, by Lemma~\ref{le5.2}, it is
sufficient to show that
$\int_{v \ge0} \log v \,d\hat{\mu}_u(v)$ is uniformly bounded in $u$,
which would follow from
the uniform integrability of $\{ |\log V_{T_u} - \log u| \}$. To this
end,\vspace*{2pt} we apply the corollary to Theorem~2 of \citet{TLDS79}.
Note that $V_{T_u} = \widetilde{V}_{T_u}$, where $\widetilde{V}_n = S_n + \delta_n$
for a sequence $\{ \delta_n \}$ which is slowly changing [cf.
\citet{JCAV11}, Lemma 4.1].
Also, using \citet{JCAV11}, Lemma 5.5, it is easy to verify that
%
%e5.33 #&#
\begin{equation}
\label{july1c} \frac{\xi}{2} {\mathbf E}_\xi\bigl[ \llvert
\delta_{T_u}-\delta_{T_u-1}\rrvert{\mathbf1}_{\{T_u < K \}}
\bigr] \le{\mathbf E}_\xi\bigl[ \log\bigl(\widebar{Z}^{(p)}
\bigr)^\xi\bigr] < \infty.
\end{equation}
Note that conditions (6)--(8) of \citet{TLDS79} are also satisfied
with $\alpha=1$. In this regard, notice that Theorem~2 of their article
is actually valid if their
equation (8) is replaced by uniform continuity in probability of $\{
\delta_n \}$, as given in equation (4.2) of \citet{MW82},
and the latter condition holds since
$\delta_n$ converges w.p.1 to a proper r.v.
We conclude
$\{ |\log V_{T_u} - \log u| \}$ is uniformly integrable.
Then~\eqref{july1b} follows since
$\hat{\mu}_u \Rightarrow\hat{\mu}$.

Finally, applying the dominated convergence theorem to the second term
in~\eqref{july1a}
and invoking~\eqref{lm4230},
we conclude
$({\log u})^{-1} {\mathbf E} [ L(M) | (V_0/u) \sim\hat{\mu} ]
\to1/|\Lambda^\prime(0)|$, as required.
\end{pf*}

%
% PROOF OF OPTIMALITY
%

%s6 #&#
\section{Proof of optimality}\label{sec6} The idea of the proof is similar to
\citet{JC02}, Theorem~3.4, but new technical issues arise since we deal with a process generated
by~\eqref{intro1} rather than a random
walk process.

%
% PROOF OF THEOREM 2.4
%

\begin{pf*}{Proof of Theorem~\ref{THoptimality}}
Let $\bfnu\in{\mathfrak M}$. First we show that
%
%e6.1 #&#
\begin{equation}
\label{newTh231} \liminf_{u \to\infty} \frac{1}{\log u} {\mathbf
E}_{\bfnu} \bigl[ \bigl({\mathcal E}_{u}^{(\bfnu)}
\bigr)^2 \bigr] \ge-2\xi.
\end{equation}

To establish~\eqref{newTh231}, set
\[
\mu_{\mathfrak D}(E;w,q) = \cases{ \mu_\xi(E), &\quad$E \in{
\mathcal B}\bigl({\mathbb R}^3\bigr)$, $w \in{\mathbb R}$ \mbox{ and
} $q= 0$;
\vspace*{2pt}\cr
\mu(E), &\quad$E \in{\mathcal B}\bigl({\mathbb R}^3
\bigr)$, $w \in{\mathbb R}$ \mbox{ and } $q= 1$.}
\]
(Intuitively, $w$ corresponds to the level of the process $\{ \log
V_{n-1}/\log u \}$,
while \mbox{$q=1$} indicates that
$\{ V_n \}$ has exceeded level $u$ by the previous time.)

If $\nu\ll\mu_{\mathfrak D}$, then by a standard argument [cf.
\citet{JC02}, equations (4.54), (4.55)], utilizing the
Radon--Nikodym theorem,
\begin{eqnarray*}
{\mathbf E}_{\bfnu} \bigl[ \bigl({\mathcal E}_{u}^{(\bfnu)}
\bigr)^2 \bigr]&:=& {\mathbf E}_\bfnu\Biggl[
N_u^2 {\mathbf1}_{\{T_u < K\}} \prod
_{i=1}^{K} \biggl( \frac{d\mu}{d\nu}(Y_i;W_i,Q_i)
\biggr)^2 \Biggr]
\\
& =& {\mathbf E}_{\mathfrak D} \Biggl[ N_u^2 {\mathbf
1}_{\{T_u < K\}} \prod_{i=1}^{K} \biggl(
\frac{d\mu}{d\mu_{\mathfrak D}}(Y_i; W_i, Q_i)
\biggr)^2 \frac{d\mu_{\mathfrak D}}{d\nu} (Y_i; W_i,Q_i)
\Biggr].
\end{eqnarray*}
Note $\frac{d\mu}{d\mu_{\mathfrak D}} = \frac{d\mu}{d\mu_\xi}$ for $Q_i
= 0$,
while
$\frac{d\mu}{d\mu_{\mathfrak D}} = 1$ for $Q_i = 1$.
Hence
%
%e6.2 #&#
\begin{equation}
\label{pfthm233}\qquad {\mathbf E}_{\bfnu} \bigl[ \bigl({\mathcal
E}_{u}^{(\bfnu)} \bigr)^2 \bigr] = {\mathbf
E}_{\mathfrak D} \Biggl[ N_u^2 {\mathbf
1}_{\{T_u < K\}} \prod_{i=1}^{T_u} \biggl(
\frac{d\mu}{d\mu_\xi}(Y_i) \biggr)^2 \prod
_{j=1}^K \frac{d\mu_{\mathfrak D}}{d\nu} (Y_j;
W_j,Q_j) \Biggr].
\end{equation}
Thus setting
\[
U_i = \log\biggl( \frac{d\nu}{d\mu_{\mathfrak D}}(Y_i;
W_i,Q_i) \biggr) \quad\mbox{and} \quad R_n = \sum_{i=1}^n
U_i,
\]
we conclude by Jensen's inequality that
%
%e6.3 #&#
\begin{eqnarray}
\label{pfthm234} {\mathbf E}_{\bfnu} \bigl[ \bigl({\mathcal
E}_{u}^{(\bfnu)} \bigr)^2 \bigr] &=& {\mathbf
E}_{\mathfrak D} \bigl[ N_u^2 {\mathbf1}_{\{ T_u < K \}}
e^{-2\xi S_{T_u} - R_{K}} \bigr]
\nonumber\\[-8pt]\\[-8pt]
& \ge& p_u \exp\bigl\{ {\mathbf E}_{\mathfrak D} [ -2\xi
S_{T_u} - R_{K} | T_u < K ] \bigr\},\nonumber
\end{eqnarray}
where
$p_u:= {\mathbf P}_\xi \{ T_u < K \} \to\Theta> 0$ as $u \to\infty$.
It follows from~\eqref{pfthm234} that
%
%e6.4 #&#
\begin{eqnarray}
\label{pfthm236} \liminf_{u \to\infty} \frac{1}{\log u} \log{\mathbf
E}_{\bfnu} \bigl[ \bigl({\mathcal E}_{u}^{(\bfnu)}
\bigr)^2 \bigr] & \ge& -\limsup_{u \to\infty}
\frac{1}{\log u} {\mathbf E}_\xi[ 2 \xi S_{T_u} {\mathbf
1}_{\{T_u < K\}} ]
\nonumber\\[-8pt]\\[-8pt]
&&{} - \limsup_{u \to\infty}
\frac{1}{\log u} {\mathbf E}_{\mathfrak D} [ R_K {\mathbf
1}_{\{ T_u < K\}} ].\nonumber
\end{eqnarray}

To identify the first term on the RHS of~\eqref{pfthm236}, note by
Wald's identity that
%
%e6.5 #&#
\begin{equation}
\label{march27a} {\mathbf E}_\xi[ \log A ] {\mathbf E}_\xi[
T_u \wedge K ] = {\mathbf E}_\xi[ S_{T_u} {\mathbf
1}_{\{T_u < K \}} ] + {\mathbf E}_\xi[ S_{K} {\mathbf
1}_{\{ K \le T_u \}} ].
\end{equation}
Now\vspace*{2pt} $(\log u)^{-1} {\mathbf E}_\xi[ T_u \wedge K] \to(\log u)^{-1} {\mathbf
E}_\xi[T_u | T_u < K]$ as $u \to\infty$
(by Theorem~\ref{MainTh-2}). Also,
$ {\mathbf E}_\xi [ S_{K} {\mathbf1}_{\{ K \le T_u \}} ] \to{\mathbf E} [ S_K
e^{-\xi S_K} {\mathbf1}_{\{ K < \infty\}} ]$ as $u \to\infty$,
which is obviously finite on $\{ S_K > 0 \}$, and which is finite on
$\{ S_K \le0 \}$ since\break  (as $e^{\xi x}\geq 1 + \xi x, x>0$) it can be bounded
by a constant multiple of\break
${\mathbf E}_\xi [ e^{-\xi S_K} {\mathbf1}_{\{S_K \le0, K < \infty\}} ] =
{\mathbf E} [ {\mathbf1}_{\{S_K \le0, K < \infty\}} ] < \infty$.
Thus,\vspace*{1pt} using that ${\mathbf E}_\xi [ \log A ] = \Lambda^\prime(\xi)$, it
follows from Theorem~\ref{MainTh-2}, equation~\eqref{thm23b}, and the above
discussion that the middle term
of~\eqref{march27a}
must satisfy
%
%e6.6 #&#
\begin{equation}
\label{pfthm2321} \lim_{u \to\infty} \frac{1}{\log u} {\mathbf
E}_\xi[ S_{T_u} {\mathbf1}_{\{ T_u < K \}} ] = 1.
\end{equation}

To handle the second limit on the RHS of~\eqref{pfthm236}, first assume,
for the moment, that $\log( \frac{d\nu}{d\mu_{\mathfrak D}} )$
is bounded from below by a finite constant. This assumption will later
be removed.
Recall that $U_i = \log ( \frac{d\nu}{d\mu_{\mathfrak D}}(Y_i;W_i,Q_i)
)$ and $R_n = \sum_{i=1}^n U_i$.\vspace*{1pt}
Now it follows by an application of Jensen's inequality that
%
%e6.7 #&#
\begin{eqnarray}
\label{jensen1} {\mathbf E}_{\mathfrak D} \bigl[ U_n |
(W_{n}, Q_{n})=(w,q) \bigr] &=& \int_{{\mathbb R}^3}
\log\biggl( \frac{d\nu}{d\mu_{\mathfrak D}} (y;w,q) \biggr) \,d\mu
_{\mathfrak D}(y;w,q)
\nonumber\\[-8pt]\\[-8pt]
&\le&\log\int_{{\mathbb R}^3} \,d\nu(y;w,q) = 0\nonumber
\end{eqnarray}
[where we have suppressed the dependence on $(w,q)$ in the above
integrals], and
consequently, after a short argument, we conclude that
\[
{\mathcal M}_n:= R_n {\mathbf E}_{\mathfrak D} [ {\mathbf
1}_{\{ T_u < K \}} | {\mathfrak F}_n ]
\]
is a supermartingale.
Hence by the optional sampling theorem,
%
%e6.8 #&#
\begin{equation}
\label{march27b} \limsup_{u \to\infty} \frac{1}{\log u} {\mathbf
E}_{\mathfrak D} [ R_K {\mathbf1}_{\{T_u < K \}} ] \le0.
\end{equation}
Then~\eqref{newTh231}
follows from~\eqref{pfthm2321} and~\eqref{march27b}.
If $\log ( \frac{d\nu}{d\mu_{\mathfrak D}} )$ is not bounded from
below by a constant, then we can replace $\nu$ with a larger
measure, $\nu^{(\varepsilon)}:= \nu+ \varepsilon\mu_{\mathfrak D}$, where
$\varepsilon> 0$. Then the entire proof can be repeated without
significant change, and we again conclude~\eqref{newTh231} upon letting
$\varepsilon\downarrow0$. We omit the details, which are straightforward.

Next, we show that strict inequality holds in~\eqref{newTh231} when
$\bfnu\in{\mathfrak M}$ differs from the dual measure.
Now if $\nu\neq\mu_{\mathfrak D}$, then, in view of
\eqref{jensen1},
there exists a point $(w,q)$ where
%
%e6.9 #&#
\begin{equation}
\label{pfthm2341} {\mathbf E}_{\mathfrak D} [ U_n | W_{n} =
w, Q_{n} = q ] = -2\Delta\qquad\mbox{for some } \Delta> 0.
\end{equation}
Then, from the definition of $U$ and an application of the
Radon--Nikodym theorem, it follows from the continuity assumption
(C$_0$) that
for some neighborhood $G$ of~$w$,
%
%e6.10 #&#
\begin{equation}
\label{march28a} {\mathbf E}_{\mathfrak D} [ U_n | W_n = w,
Q_n = q ] \le- \Delta,\qquad w \in G.
\end{equation}

We now show that by sharpening the estimate in Jensen's inequality on
the set $G \times\{ q\}$, we obtain a \textit{strict} inequality in \eqref
{newTh231}.
As before, we begin by assuming that $\log ( \frac{d\nu}{d\mu
_{\mathfrak D}} )$ is bounded from below by a constant.
Then by repeating our previous argument, but using the sharper estimate
\eqref{march28a} when $w \in G$ and $q$ given as in~\eqref{pfthm2341},
together with Jensen's
inequality for the remaining values of $(w,q)$, we obtain that
\begin{eqnarray*}
{\mathcal M}_n^\ast&:=& \bigl( U_1^\ast
+ \cdots+ U_n^\ast\bigr) {\mathbf E}_{\mathfrak D} [{\mathbf
1}_{ \{T_u < K \}} | {\mathfrak F}_n ],
\\
U_i^\ast&:=& U_i + \Delta{\mathbf
1}_{\{W_n \in G\}} {\mathbf1}_{\{ Q_{n} =q\}},
\end{eqnarray*}
is a supermartingale.
Applying the optional sampling theorem, we deduce that
%
%e6.11 #&#
\begin{equation}
\label{pfthm2342} {\mathbf E}_{\mathfrak D} [ R_K {\mathbf
1}_{\{ T_u < K \}} ] \le- \Delta\bigl\{ {\mathbf1}_{\{q=0\}} {\mathbf
E}_{\mathfrak D} \bigl[ {\mathcal O}_u^{(0)} \bigr] + {
\mathbf1}_{\{q=1\}} {\mathbf E}_{\mathfrak D} \bigl[ {\mathcal
O}_u^{(1)} \bigr] \bigr\},
\end{equation}
where
\[
{\mathcal O}_u^{(0)}:= \sum
_{n=0}^{T_u} {\mathbf1}_{\{ W_n \in G \}}
\]
and
\[
{\mathcal O}_u^{(1)}:= \sum
_{n=T_u+1}^{K} {\mathbf1}_{\{W_n \in G \}}.\vadjust{\goodbreak}
\]
Note that ${\mathcal O}_u^{(0)}$ denotes the occupation time which the
scaled process\break  $\{\log V_{n}/\log u \}$ spends in the interval $G$
during a trajectory starting at time 0 and ending at time $T_u$, while
${\mathcal O}_u^{(1)}$ denotes the occupation time that $\{\log
V_{n}/\log u \}$ spends in the interval $G$ during a trajectory
starting at time $T_u$ and ending at time $K$. Note that for all $n \in
\pintegers$,
\[
V_{0} - W \le\frac{V_n}{A_1 \cdots A_n} \le V_0 + W\qquad\mbox{where } W:= \sum_{i=1}^\infty
\frac{|B_i| + A_i |D_i|}{A_1 \cdots A_i}.
\]
Now suppose that $G^\prime:= [u^{s^\prime},u^{t^\prime}] \subset
[u^s,u^t] \subset G$, where $s
< s^\prime< t^\prime< t$. Then in the $\xi$-shifted measure,
the transient process $\{V_n \}$ enters $G^\prime$ w.p. $p_u \to\Theta
>0$. Now, in the previous equation, take $V_0$ to be the position of
this process at its first passage time into $G^{\prime}$, so that $V_0
\ge u^{s^\prime}$.
Since $W$ is a proper r.v. w.p.1 in the $\xi$-shifted measure,
it follows that for some $\varepsilon> 0$,
${\mathbf P} \{ (V_0+W)/(V_0 - W) - 1 > u^{-\varepsilon} \} \to0$
as $u \to\infty$ (and an analogous estimate holds when $V_0 < W$).
Thus we see that $\{ \log V_n \}$ is well-approximated
by $\{ S_n\}$. Since, as a multiplicative random walk, the occupation
time of $\{ e^{S_n} \}$ in $G^\prime$
is at least $c\log u$ for some $c>0$, we conclude (after a short
argument) that
%
%e6.12 #&#
\begin{equation}
\label{march28c} \liminf_{u \to\infty} \frac{1}{\log u} {\mathbf E}
\bigl[ {\mathcal O}_u^{(0)} \bigr] \ge\eta>0.
\end{equation}
Substituting this estimate into~\eqref{pfthm2342} yields, for the case
$q=0$ in~\eqref{pfthm2341}, that
%
%e6.13 #&#
\begin{equation}
\label{pfthm2354} \limsup_{u \to\infty} \frac{1}{\log u} {\mathbf
E}_{\mathfrak D} [ R_K {\mathbf1}_{\{ T_u < K \}} ] \le-\Delta\eta
<0.
\end{equation}
Now substituting~\eqref{pfthm2354} and~\eqref{pfthm2321} into
\eqref{pfthm236},
we obtain that the LHS of~\eqref{pfthm236} is $\ge-2\xi+ \Delta\eta
$, as required.

If $q=1$ in~\eqref{pfthm2341}, the argument is similar. Here we study a
trajectory in the original measure, beginning at the level $V_{T_u}$
and returning to the set $K$. Setting $V_0 \stackrel{\mathcal D}{=}
V_{T_u}$, then we may again observe that $\{ \log V_n \}$ behaves
similarly to a random walk or, more precisely,
%
%e6.14 #&#
\begin{equation}
\label{pfthm231n} \sup_n \Biggl| V_n -
V_0 \prod_{i=1}^n
A_i \Biggr| \le W^\prime\qquad\mbox{where } W^\prime:=
\sum_{i=1}^\infty\widetilde{B}_i
\prod_{j=i+1}^\infty A_j
\end{equation}
as long as $\{ V_0,\ldots, V_{n-1} \}$ is nonnegative.
Then by a straightforward argument based on the law of large numbers,
%
%e6.15 #&#
\begin{equation}
\label{pfthm235n} \liminf_{u \to\infty} \frac{1}{\log u} {\mathbf E}
\bigl[ {\mathcal O}_u^{(1)} \bigr] \ge\tilde{\eta} > 0
\end{equation}
and so we obtain that the LHS of~\eqref{pfthm236} is $\ge-2\xi+
\Delta\tilde{\eta}$. (For more details, see our preprint under the same title in Math arXiv.)

If $\log ( \frac{d\nu}{d\mu_{\mathfrak D}} )$ is not bounded from
below by a constant, then replace $\nu$ with
$\nu^{(\varepsilon)}:= \nu+ \varepsilon\mu_{\mathfrak D}$, where $\varepsilon
> 0$, and the proof carries through with little modification.
Finally, to complete the proof of theorem, note that if we do not have
$\nu\ll\mu_{\mathfrak D}$, as we have assumed throughout this proof,
then by an application of the Radon--Nikodym theorem,
$\nu= \nu_a + \nu_s$, where $\nu_a \ll\mu_{\mathfrak D}$ and $\nu_s
\perp\mu_{\mathfrak D}$.
The proof can now be repeated, replacing everywhere $\nu$ with $\nu_a$;
cf. \citet{JC02}, proof of Theorem~3.4.
We omit the details.
\end{pf*}

% zodis "Acknowledgments" paliekamas pagal autoriu

%suskaldyti doi

% imsref loaded by linak, 2014-02-27 09:23:30
%
% imsref loaded by linak, 2014-03-26 10:45:59

\printaddresses

\end{document}